\documentclass[12pt]{article}
\usepackage{amsmath}
\usepackage{amssymb}
\usepackage{amsfonts}
\usepackage{amscd}
\usepackage{enumerate}
\usepackage{fancyheadings}
\usepackage{color}
\usepackage{graphicx}
\newtheorem{Lem}{Lemma}
\newtheorem{Th}{Theorem}
\newtheorem{Def}{Definition}
\newtheorem{Prop}{Proposition}
\newtheorem{Assume}{Assumption}
\DeclareMathOperator{\supp}{supp}
\DeclareMathOperator{\graph}{graph}

\DeclareMathOperator{\diag}{diag}
\begin{document}

\begin{center}
\bf{\Large Structure of the Semi-Classical Amplitude for General Scattering Relations}
\end{center}

\noindent
{\bf Ivana Alexandrova}

\noindent
Department of Mathematics, University of Toronto, Toronto, Ontario, Canada 
M5S 3G3, Tel.: 1-416-946-0318, Fax: 1-416-978-4107, email: alexandr@math.toronto.edu

\noindent
July 26, 2004

\begin{abstract}
We consider scattering by general compactly supported
semi-classical perturbations of the Euclidean Laplace-Beltrami operator.
We show that if the suitably cut-off resolvent of the Hamiltonian
quantizes a Lagrangian relation on the product cotangent
bundle, the scattering amplitude quantizes the natural scattering
relation. 
In the case when the resolvent is tempered, which 
is true under some non-resonance assumptions, and when we
work microlocally near a non-trapped ray our result implies that
the scattering amplitude defines a semiclassical Fourier integral
operator associated to the scattering relation in a neighborhood
of that ray. 
Compared to previous work we allow this relation to 
have more general geometric structure.
\end{abstract}
\section{Introduction and Statement of Results}
We study the semi-classical scattering amplitude at non-trapping energies
for compactly supported metric and potential perturbation
$P\left(h\right)$ of
the Euclidean Laplace-Beltrami operator $P_{0}\left(h\right)$ on
$\mathbb{R}^{n}$.
The scattering amplitude at energy $\lambda>0$ is the amplitude of the
leading term in the
asymptotic expansion of an outgoing solution of
$\left(P\left(h\right)-\lambda\right)u=0$ as
$\|x\|=r\to
\infty$.
We prove that the scattering amplitude quantizes the natural
scattering relation in a sense of global 
semi-classical Fourier Integral Operators. 

\subsection{A Survey of Earlier Results}

The structure of the scattering amplitude has been of considerable interest to researchers 
in mathematical physics.
To outline the earlier results, we begin by making some definitions.
Let $P(h)=-\frac{1}{2}h^{2}\Delta+V,$ where
\begin{equation}\label{potential}
\left|\frac{\partial^{\alpha}}{\partial x^{\alpha}}V(x)\right|\leq C_{\alpha}(1+|x|)^{-\mu-|\alpha|},\; x\in\mathbb{R}^{n}, \mu>0.
\end{equation}
Let $\omega_0\in\mathbb{S}^{n-1}$ and for $z\in\omega^{\perp}_{0},$ with $z$ its coordinate 
representation in 
$\mathbb{R}^{n-1},$
let \[\left\{q_{\infty}\left(t; z, \lambda\right), p_{\infty}\left(t; z, 
\lambda\right)\right\}\] be the 
unique phase 
trajectory such that 
\begin{gather*}
\lim_{t\to -\infty}\left|q_{\infty}\left(t; z, 
\lambda\right)-\sqrt{2\lambda}\omega_0 
t-z\right|=0,\\
\lim_{t\to -\infty}\left|p_{\infty}\left(t; z, 
\lambda\right)-\sqrt{2\lambda}\omega_0\right|=0
\end{gather*}
in the $C^{\infty}$ topology for the impact parameter $z.$
Then there exist $\xi_{\infty}\left(z; \lambda\right), r_{\infty}\left(z; 
\lambda\right) \in C^{\infty}(\omega^{\perp})$ with $\xi_{\infty}(z; \lambda)\in\mathbb{S}^{n-1},$ such that
\begin{gather*}
\lim_{t\to\infty}\left|q_{\infty}\left(t; z, 
\lambda\right)-\sqrt{2\lambda}\xi_{\infty}\left(z; 
\lambda\right)t-r_{\infty}\left(z; \lambda\right)\right|=0,\\
\lim_{t\to\infty}\left|q_{\infty}\left(t; 
z,\lambda\right)-\sqrt{2\lambda}\xi_{\infty}\left(z; 
\lambda\right)\right|=0,
\end{gather*}
in the $C^{\infty}$ topology for $z.$
The trajectory $\left\{q_{\infty}\left(t; z, \lambda\right), p_{\infty}\left(t; z, \lambda\right)\right\}$ is then said to have initial direction $\omega_0$ and final direction $\theta_0=\xi_{\infty}(z, \lambda)$ and $\theta_0$ is said to be non-degenerate, or regular, for $\omega_0$ if for all 
$z\in\omega^{\perp}_{0}$ with $\xi_{\infty}\left(z; \lambda\right)=\theta_0,$ 
the angular density $\hat{\sigma}\left(z; \lambda\right)$ for the trajectory 
$\left\{q_{\infty}\left(t; z, \lambda\right), 
p_{\infty}\left(t; z, \lambda\right)\right\}$ satisfies
\begin{equation}\label{angulard}
\hat{\sigma}\left(z; 
\lambda\right):=\left|\det\left(\xi_{\infty}, \frac{\partial}{\partial z_{1}}\xi_{\infty}, \dots, \frac{\partial}{\partial z_{n-1}}\xi_{\infty} \right) \right|\neq 0.
\end{equation} 

The first asymptotic expansion of the semi-classical scattering amplitude was given by 
Vainberg \cite{V}.
He considers the semi-classical Schr\"{o}dinger operator with a potential $V\in 
C_{c}^{\infty}(\mathbb{R}^{n}; \mathbb{R})$ and assumes \eqref{angulard}.
For the associated scattering amplitude at energies $\lambda> \sup V$ he then proves an asymptotic expansion of the scattering amplitude in the form
\begin{equation}\label{vexpansion}
f(\theta, \lambda)=\sum_{j=1}^{l}\hat{\sigma}\left(z_{j};
\lambda\right)^{-1/2}\exp\left(ih^{-1}S_{j}-i\mu_{j}\pi/2\right)+\mathcal{O}\left(h\right),
\end{equation}
where $\left(z_{j}\right)_{j=1}^{l}\equiv\left(\xi_{\infty}^{-1}(\cdot; \lambda)\right)(\theta_0),$  $S_j$ 
is a modified action along the $j-$th $(\omega_0, \theta)$ trajectory, $\theta\in U,$ where $U\subset\mathbb{S}^{n-1}$ is a sufficiently small open neighborhood of $\theta_0,$
and $\mu_{j}$ is the 
path index of the trajectory.
The error term is estimated uniformly in $\theta\in U.$
 
Majda \cite{AM} considers scattering processes defined by the classical 
wave equation in 
the presence of a convex obstacle in $\mathbb{R}^{n}$ where $n=2, 3.$
For Dirichlet, Neumann, and in the case of three-dimensional space, 
impedance boundary 
conditions, he proves an asymptotic expansion of the scattering amplitude, 
the leading term 
of which is the product of the Gauss curvature and the reflection 
coefficient evaluated at a 
point on the boundary of the obstacle.
In this setting the scattering amplitude is the coefficient of the leading 
term 
in the asymptotic 
expansion of an outgoing solution of the reduced wave equation.
To establish the aforementioned result, he studies the radiation 
pattern of a solution which approximates this 
outgoing solution 
near the boundary of the obstacle.
He also applies his main result to inverse scattering 
problems for convex bodies with the above boundary conditions.
In particular, he proves that both the shape of the boundary of the 
obstacle and the nature 
of the boundary conditions are completely determined by the asymptotic 
limit of the 
scattering amplitude.

Guillemin \cite{G} discusses similar asymptotic expansions.
He studies the behavior of the scattering matrix in several different 
settings: on a compact manifold, in obstacle scattering, for a compactly 
supported perturbation of the Euclidean metric on $\mathbb{R}^{n},$ and 
for a 
quotient Riemannian manifold.
In each case, he presents formulas for the kernel of the scattering matrix 
at energy $\lambda$ of the form
\begin{equation*}
S\left(\lambda, \omega, \theta\right)=c\left(\lambda, 
n\right)\sum_{j=1}^{N} c\left(j\right) 
\left|J_{j}(\lambda)\right|^{-1/2}e^{i\lambda 
T_{j}}+\mathcal{O}\left(\frac{1}{\lambda}\right), \theta\ne\omega
\end{equation*} 
under the assumption that there are $N$ scattering rays with initial 
direction $\omega$ and final direction $\theta,$ where $T_{j}$ is the 
sojourn time of the $j-$th scattering ray and $J_{j}$ is the scattering 
differential cross-section evaluated at the point of incidence of the 
$j-$th scattering ray.
In the case of the quotient, the scattering matrix is a unitary matrix 
of size depending on the topology of the manifold. 
For each energy level $\lambda$ its $jk-$th entry has the form
\begin{equation*}
S_{jk}\left(\lambda\right)=ac\left(\lambda\right)\sum 
e^{-T_{k}\left(i\lambda+1/2\right)},
\end{equation*}
where
\begin{equation*}
c\left(\lambda\right)=\int^{\infty}_{-\infty}\frac{dq}{\left(1+q^{2}\right)^{1/2+i\lambda}}.
\end{equation*}

To derive these results, Guillemin uses the representation of the 
scattering operator in terms of the wave operators.
He also derives a formula for the scattering differential cross-section in 
the case of scattering by a smooth convex obstacle from which he 
deduces that the asymptotic behavior of the scattering amplitude 
determines 
the shape of the scatterer.

A different form of the asymptotic expansion of the scattering matrix was given by Protas 
\cite{P}.
He works in the setting considered by Vainberg \cite{V} and proves an 
asymptotic expansion of the scattering amplitude in terms of canonical Maslov operators.
This expansion holds for a fixed initial direction and uniformly in 
an open set 
containing the final direction and disjoint from the initial direction.

Yajima \cite{Y} was the first to prove an asymptotic expansion of the form \eqref{vexpansion} of the scattering amplitude
for potential perturbations $V\in C^{\infty}(\mathbb{R}^{n}; \mathbb{R})$ of the semi-classical Laplacian satisfying \eqref{potential} for a constant $\mu>\max\left(1, \frac{n-1}{2}\right).$ 
He also works at non-trapping energies and with outgoing directions non-degenerate for the initial directions.
His results, however, are only valid in the $L^{2}$ sense and under the non-trapping 
assumption and for outgoing directions which are non-degenerate for the fixed incoming direction.

Robert and Tamura \cite{RT} work in the same setting as Yajima \cite{Y} with $\mu\geq 1.$
For scattering amplitude at non-trapping energies $\lambda>0,$ which now satisfies 
\[f\left(\lambda, h\right)\in 
C^{\infty}\left(\mathbb{S}^{n-1}\times\mathbb{S}^{n-1}\backslash\diag
\left(\mathbb{S}^{n-1}\times\mathbb{S}^{n-1}\right)\right),\] they establish an asymptotic expansion of the form \eqref{vexpansion} with
\begin{equation}\label{phaseRT}
S_{j}=\int_{-\infty}^{\infty}\left(|p_{\infty}\left(t; z_{j}, 
\lambda\right)|^{2}/2-V\left(q_{\infty}\left(t; z_{j}, 
\lambda\right)\right)-\lambda\right) 
dt-\left\langle r_{\infty}\left(w_{j}; 
\lambda\right),\sqrt{2\lambda}\theta\right\rangle.
\end{equation}

Michel \cite{Michel} works in the same setting as Robert and Tamura but he 
allows the energy level to be trapping while satisfying the condition
\begin{equation*}
\begin{aligned}
& \text{There exists a neighborhood } W \text{ of } 
\omega_0\in\mathbb{S}^{n-1} \text{ such}\\
&\text{ that } \forall \omega'\in W, \forall 
z\in(\omega')^{\perp}\lim_{t\to\infty}|q_{\infty}\left(t, z, 
\omega'\right)|=\infty.
\end{aligned}
\end{equation*}
He further assumes that there exists $\epsilon>0$ such that the resonances 
$\lambda_{j}$ satisfy $|\Im \lambda_{j}|\geq Ch^{q}$ for $\Re 
\lambda_{j}\in\left[\lambda-\epsilon, \lambda+\epsilon\right].$
Under these assumptions he establishes the same asymptotic expansion of 
the scattering amplitude (\ref{vexpansion}).
Like Robert and Tamura, Michel also uses Isozaki-Kitada's representation 
formula of the scattering amplitude.

\subsection{Statement of Main Theorem}
In this article we analyze the semi-classical scattering amplitude
from a different point of view and {\em without} making any geometric
assumptions on the scattering relation such as (\ref{angulard}).
We show that under a microlocal assumption on the resolvent, essentially 
the assumption that a suitably cut-off resolvent quantizes the flow 
relation, the scattering amplitude is a semi-classical
Fourier Integral Operator associated to the classical scattering
relation. 
In other words, it quantizes that canonical relation.

We will work in the abstract ``black box'' framework of \cite{SZ}
which means that we can formulate our hypotheses independently of
the structure of the scatterer. 
We refer to Section \ref{sreprsa} for the definition of the scattering amplitude $A(\lambda, h),$ 
and to \cite[Section 3.2]{AI} for a complete characterization of the class of semi-classical Fourier 
integral distributions $I_{h}.$
We also review briefly the definition of semi-classical Fourier integral distributions in Section \ref{scanal}.
The notion of microlocal localization is also reviewed in Section \ref{scanal}.

To state our main theorem, we also let $$\pi_{1}: T^{*}\mathbb{S}^{n-1}\times
T^{*}\mathbb{S}^{n-1}\times
T^{*}\mathbb{R}^{n}\times T^{*}\mathbb{R}^{n}\to T^{*}\mathbb{S}^{n-1}\times
T^{*}\mathbb{S}^{n-1}$$ and $$\pi_{2}: T^{*}\mathbb{S}^{n-1}\times T^{*}\mathbb{S}^{n-1}\times
T^{*}\mathbb{R}^{n}\times T^{*}\mathbb{R}^{n}\to T^{*}\mathbb{R}^{n}\times
T^{*}\mathbb{R}^{n}$$ denote the canonical projections.
We further introduce the following Lagrangian submanifold of $$T^*
\left(\mathbb{S}^{n-1} \times \mathbb{S}^{n-1} \times \mathbb{R}^{n} \times \mathbb{R}^{n}\right)$$
depending
on the real energy $ \lambda>0:$
\begin{gather*}
 \Lambda ( \lambda ) = \left\{ (m , d f (m)) \; : \; m \in
\mathbb{S}^{n-1} \times \mathbb{S}^{n-1} \times \mathbb{R}^{n}\backslash B(0, R_0) \times \mathbb{R}^{n}\backslash B(0, R_0) \right\} \\
f ( m ) = - \sqrt{2\lambda}\langle x , \omega \rangle+\sqrt{2\lambda}\langle y , \theta \rangle,
\ \ m = ( \omega, \theta, x, y ).
\end{gather*}
Lastly, we make the following two assumptions:
\begin{Assume}\label{assumer}
There exists $s\in\mathbb{R}$ such that for every $\varphi\in C_{c}^{\infty}(\mathbb{R}^{n}),$
$\varphi\equiv \text{ const. }$ on $B(0, R_{0}),$ $\left\|\varphi R(\lambda,
h)\varphi\right\|_{\mathcal{B}(L^{2}(\mathbb{R}^{n}))}=\mathcal{O}(h^{s}).$
\end{Assume}
\begin{Assume}\label{assumerfio}
There exists a Lagrangian submanifold $\Lambda_R(\lambda)\subset 
\pi_{2}\left(\Lambda(\lambda)\right)$ of $ T^* {\mathbb R}^n \times T^* 
{\mathbb R}^n $ such that 
\begin{gather*}
\text{for every }\tilde{\chi}_j\in 
C_{c}^{\infty}(\mathbb{R}^{n}\backslash
B(0, R_{0})), j=1, 2, 
\supp\tilde{\chi}_{1}\cap\supp\tilde{\chi}_{2}=\emptyset,\\
DK_{\tilde{\chi}_2 R\left(\lambda, h\right)\tilde{\chi}_1}\in 
I_{h}^{r}\left(\mathbb{R}^{n}\times
\mathbb{R}^{n}, \Lambda_R ( \lambda ) \right),
\end{gather*} 
where $D\in\Psi_{h}^{0}(1, \mathbb{R}^{n}\times\mathbb{R}^{n})$ is a microlocal cut-off to a neighborhood of 
\[
\left\{\left( x_0+t \omega_{0}, y_0+ s \theta_{0},  \sqrt{2 \lambda} \omega_{0}, -\sqrt{2 \lambda} \theta_{0}\right): s\in [s_1, s_2], t\in [t_1, t_2]\right\},\]
for some fixed $s_1<s_2<0,$ $t_1>t_2>0,$ $\omega_{0},$ $\theta_{0},$ $x_0,$ and $y_0,$ where $K_{\tilde{\chi}_2 R\left(\lambda, h\right)\tilde{\chi}_1}$ denotes the Schwartz kernel of the cut-off resolvent.
\end{Assume}
The last assumption means that the cut-off resolvent is a Fourier 
Integral Operator microlocally near some incoming and outgoing directions.
The first assumption is made so that the notion of applying semi-classical 
pseudodifferential operators (which here are defined up to residual terms
in ${\mathcal O} ( h^\infty ) $) makes sense. It will be used 
explicitly in the discussion of the resolvent near a non-trapped
trajectory in Sections \ref{snt} and \ref{str}.

We also note that implicit in our assumptions
is the fact that $ \lambda $ is non-resonant, in the sense that
the resolvent does not have a pole at $ \lambda.$

We can now state our

\medskip
\noindent
{\bf Main Theorem.}
{\em Suppose that Assumptions 1 and 2 hold. 

Then $$\pi_1 \circ
\left(\pi_2 |_{\Lambda( \lambda ) } \right)^{-1}
\Lambda_R ( \lambda )$$ is a smooth Lagrangian submanifold of
$ T^* S^{n-1} \times T^* S^{n-1} $ near 
\begin{equation*}
\bar{p}=\left(\omega_{0}, -\sqrt{2\lambda} d_\omega \langle x_0, \omega_{0} \rangle;\;
\theta_{0}, \sqrt{2\lambda}d_{\theta}\langle y_0, \theta_{0} \rangle\right)
\end{equation*}
and for every $C\in\Psi_{h}^{0}(1, \mathbb{S}^{n-1}\times\mathbb{S}^{n-1})$ with 
symbol supported near $\bar{p}$ we have
\[CK_{A(\lambda, h)}\in I_{h}^{r+\frac{1}{2}}\left({\mathbb 
S}^{n-1} \times
{\mathbb S}^{n-1} ,
\pi_1 \circ \left(\pi_2|_{\Lambda(\lambda))}\right)^{-1} ( \Lambda_R ( \lambda))\right),\]
where $K_{A(\lambda, h)}$ denotes the Schwartz kernel of the scattering amplitude.}

\begin{figure}[t]
\begin{center}
\input{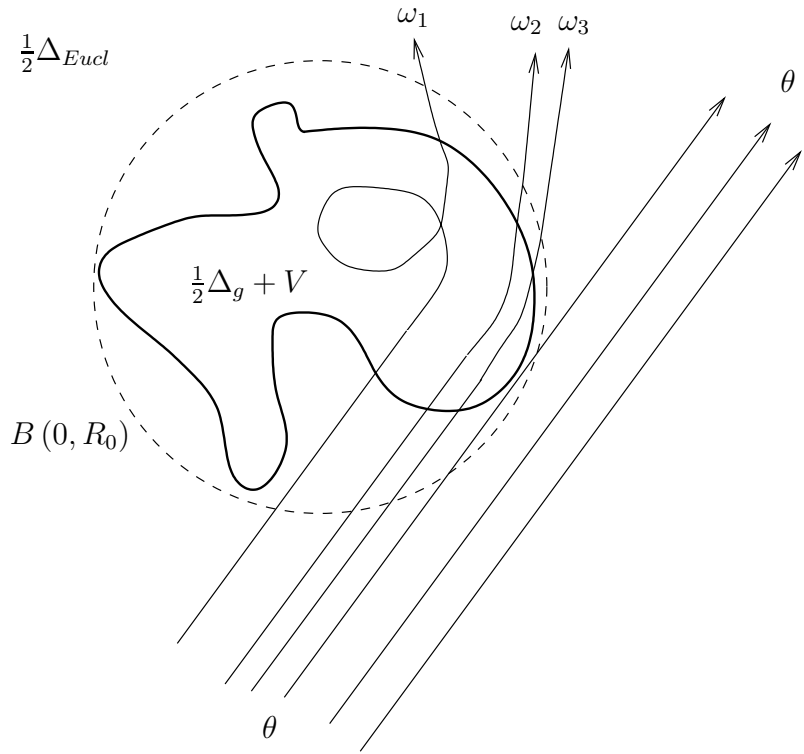}
\end{center}
\caption{A typical trajectory in the presence of a perturbation.
At the boundary between the perturbed trajectories and free trajectories
the scattering relation has degenerate projections to the $ (\theta , 
\omega ) $ variables.}
\label{fig:pert}
\end{figure}

In the special case when the non-degeneracy assumption 
\eqref{angulard} holds we recover the phases \eqref{phaseRT} in \eqref{vexpansion} -- see 
Theorem 
\ref{tmicrol} 
below. We expect that a finer analysis based on our method
would give a precise description of amplitudes. 
What is different here is the fact that we can handle the cases in which 
the scattering relation cannot be parameterized simply. That always
occurs at the transition between the perturbation and free
propagation -- see Figure~\ref{fig:pert}. 

This paper is organized as follows.
In Section \ref{snotation} we introduce some of the notation, which we 
will use throughout this article.
We review the relevant part of semi-classical analysis in Section 
\ref{scanal}.
The representation of the scattering matrix, which we will use here, is given 
in Section \ref{sreprsa}.
Section \ref{cgeom} is dedicated to the geometric aspects of the problem with the 
scattering relation defined and studied in Section \ref{scr}, and the 
resolvent relation, i.~e., the canonical relation which we will prove is 
quantized by the cut-off resolvent, described in Section \ref{sresolrel}.
The proof of the main theorem is given in Section \ref{cprm}.
In Section \ref{cappl} we discuss applications of our main theorem to 
non-trapping (Section \ref{snt}) and trapping (Section \ref{str}) smooth
compactly supported perturbations of the Euclidean Laplace-Beltrami 
operator.
For that, we prove, in Section \ref{sresolv}, that the cut-off resolvent 
for such perturbations satisfies assumption \ref{assumerfio}.
The microlocal representation of the scattering amplitude 
analogous to \eqref{vexpansion} under 
the 
non-degeneracy assumption \eqref{angulard} on the angular density is given 
in Section 
\ref{smicrol}.
Our results are applied to an inverse problem in 
Section 
\ref{sinv}.

\subsection{Notation}\label{snotation}
In this section we introduce some of the notation which we will use below.
We shall denote the Euclidean norm on $\mathbb{R}^{n},$ by $\| \cdot \|$ 
and we set $B\left(0, r\right)=\left\{x\in\mathbb{R}^{n}|\|x\|\leq 
r\right\}$ for 
$r>0.$
On any smooth manifold $M$ we denote by $\sigma$ the canonical symplectic 
form on $T^{*}M$ and everywhere below we work with the canonical 
symplectic structure on $T^{*}M.$
The canonical symplectic coordinates on $T^{*}\mathbb{R}^{n}$ will be 
denoted by $(x, \xi).$
For a function $f\in C^{\infty}\left(T^{*}M\right)$ we shall denote by 
$H_{f}$ its Hamiltonian vector field.
The integral curve of $H_{f}$ with initial conditions $(x_0, \xi_0)$ will 
be denoted by 
$\gamma(\cdot; x_0, \xi_0)=(x\left(\cdot; x_0, \xi_0\right), \xi(\cdot; 
x_0, \xi_0)).$
If $C\subset T^{*}M_1\times T^{*}M_2,$ where $M_j,$ $j=1, 2,$ are smooth 
manifolds we will use the notation $C'=\{(x, \xi; y, 
-\eta): (x, \xi; y, \eta)\in C\}.$
The Euclidean norm on $\mathbb{R}^{n}$ will be denoted by $\|\cdot\|$ and 
we set $B(a, R)=\{x\in\mathbb{R}^{n}| \|a-x\|<R\},$ for $R>0.$
For a sequentially continuous operator 
$T:C^{\infty}_{c}(\mathbb{R}^{m})\to\mathcal{D}'(\mathbb{R}^{n})$ we shall 
denote by $K_{T}$ its Schwartz kernel.
For such an operator $T$ we use $T^{t}$ to denote the operator with
Schwartz kernel $K_{T^{t}}(x, y)=K_{T}(y, x).$
Unless otherwise specified, we will use $\langle\cdot, \cdot\rangle$ to 
denote the standard 
inner product on $\mathbb{R}^{n},$ $\mathbb{C}^{n},$ and 
$L^{2}(\mathbb{R}^{n}),$ and $C$ to denote a 
positive constant, which will be allowed to change from line to line. 

\section{Preliminaries}
In this section we present some of the preliminary results we shall use 
throughout this work.
\subsection{Elements of Semi-Classical Analysis}\label{scanal}
In this section we recall some of the elements of semi-classical analysis 
which we 
will use here.
First we define two classes of symbols
\begin{equation*}
S_{2n}^{m}\left(1\right)= \left\{ a\in
C^{\infty}\left(\mathbb{R}^{2n}\times(0, h_0]\right): \forall
\alpha, \beta\in\mathbb{N}^{n}, \sup_{(x, \xi, 
h)\in\mathbb{R}^{2n}\times (0, 
h_{0}]}h^{m}\left|\partial^{\alpha}_{x}\partial^{\beta}_{\xi}a\left(x, 
\xi;
h\right)\right|\leq
C_{\alpha, \beta}\right\}
\end{equation*}
and
\begin{equation*}
S^{m, k}\left(T^{*}\mathbb{R}^{n}\right)=\left\{a\in
C^{\infty}\left(T^{*}\mathbb{R}^{n}\times(0, h_0]\right): \forall \alpha, 
\beta\in\mathbb{N}^{n}, \left|\partial^{\alpha}_{x}\partial^{\beta}_{\xi} 
a\left(x, \xi;
h\right)\right|\leq
C_{\alpha,
\beta}h^{-m}\left\langle\xi\right\rangle^{k-|\beta|}\right\},
\end{equation*}
where $h_0\in(0,1]$ and $m, k\in\mathbb{R}.$
For $a\in S_{2n}\left(1\right)$ or $a\in S^{m, k}\left(T^{*}\mathbb{R}^{n}\right)$ we define the
corresponding semi-classical pseudodifferential operator of class $\Psi_{h, t}^{m}(1, \mathbb{R}^{n})$ or $\Psi_{h, t}^{m, k}(\mathbb{R}^{n}),$ respectively, by setting
\begin{equation*}
Op_{h}\left(a\right)u\left(x\right)=\frac{1}{\left(2\pi
h\right)^{n}}\int\int e^{\frac{i\left\langle x-y,
\xi\right\rangle}{h}}a\left(x, \xi; h\right)u\left(y\right) dy d\xi, \;u\in
\mathcal{S}\left(\mathbb{R}^{n}\right),
\end{equation*}
for $t\in[0, 1]$
and extending the definition to $\mathcal{S}'\left(\mathbb{R}^{n}\right)$ by
duality
(see \cite{DS}).
Below we shall work only with symbols which admit asymptotic expansions in $h$ and with
pseudodifferential operators which are quantizations of such symbols.
For $A\in\Psi_{h, t}^{k}(1, \mathbb{R}^{n})$ or $A\in\Psi_{h, t}^{m, k}(\mathbb{R}^{n}),$ we shall use $\sigma_{0}(A)$ 
and $\sigma(A)$ to denote its principal symbol and its complete symbol, respectively.
A semi-classical pseudodifferential operator will be called of principal 
type if its 
principal symbol $a_0$ satisfies
\begin{equation}\label{prtype}
a_0=0\implies da_0\ne 0.
\end{equation}

We also define the class of semi-classical distributions 
$\mathcal{D}_{h}'(\mathbb{R}^{n})$ with which we will work here
\begin{equation*}
\begin{aligned}
\mathcal{D}'_{h}(\mathbb{R}^{n}) = & \big\{u\in C^{\infty}_{h}\left((0, 
1];
\mathcal{D}'\left(\mathbb{R}^{n}\right)\right): \forall\chi\in 
C_{c}^{\infty}\left(\mathbb{R}^{n}\right) \exists\: N\in\mathbb{N}\text{ 
and 
} C_{N}>0:\\
& \quad |\mathcal{F}_{h}\left(\chi u\right)\left(\xi, h\right)|\leq
C_{N}h^{-N}\langle\xi\rangle^{N}\big\}
\end{aligned}
\end{equation*}
where
\begin{equation*}
\mathcal{F}_{h}\left(u\right)\left(\xi,
h\right)=\int_{\mathbb{R}^{n}}e^{-\frac{i}{h}\left\langle x,
\xi\right\rangle}u\left(x, h\right)dx
\end{equation*}
with the obvious extension of this definition to 
$\mathcal{E}_{h}'(\mathbb{R}^{n}).$
We shall work with the $L^{2}-$ based semi-classical Sobolev spaces 
$H^{s}(\mathbb{R}^{n}),$ $s\in\mathbb{R},$ which consist of the distributions 
$u\in\mathcal{D}_{h}'(\mathbb{R}^{n})$ such 
that 
$\|u\|_{H^{s}(\mathbb{R}^{n})}^{2}=\frac{1}{(2\pi 
h)^{n}}\int_{\mathbb{R}^{n}}(1+\|\xi\|^{2})^{s}\left|\mathcal{F}_{h}(u)(\xi, 
h)\right|^{2}d\xi<\infty.$

We shall say that $u=v$ {\it microlocally} near an open set
$U\subset T^{*}\mathbb{R}^{n}$, if
$P(u-v)=\mathcal{O}\left(h^{\infty}\right)$ in
$C_{c}^{\infty}\left(\mathbb{R}^{n}\right)$ for
every $P\in \Psi^{0}_{h}\left(1, \mathbb{R}^{n}\right)$ such that
\begin{equation}\label{P}
WF_{h}\left(P\right)\subset \tilde{U}, \bar{U}\Subset \tilde{U}\Subset
T^{*}\mathbb{R}^{n}, \tilde{U} \text{ open}.
\end{equation}
We shall also say that $u$ satisfies a property $\mathcal{P}$  {\it
microlocally} near an open set $U\subset T^{*}{\mathbb{R}^{n}}$ if there
exists $v\in\mathcal{D}_{h}'\left(\mathbb{R}^{n}\right)$ such that $u=v$
microlocally
near $U$ and $v$ satisfies property $\mathcal{P}$.

For open sets $U, V\subset T^{*}\mathbb{R}^{n},$ the operators $T,
T'\in\Psi^{m}_{h}\left(\mathbb{R}^{n}\right)$ are said to be {\it
microlocally
equivalent} near $V\times U$ if for any $A, 
B\in\Psi_{h}^{0}\left(\mathbb{R}^{n}\right)$
such that
\begin{equation*}
WF_{h}\left(A\right)\subset\tilde{V},
WF_{h}\left(B\right)\subset\tilde{U},
\bar{V}\Subset\tilde{V}\Subset T^{*}\mathbb{R}^{n},
\bar{U}\Subset\tilde{U}\Subset T^{*}\mathbb{R}^{n}, \tilde{U}, \tilde{V}
\text{ open }
\end{equation*}
\begin{equation*}
A\left(T-T'\right)B=\mathcal{O}\left(h^{\infty}\right)\colon\mathcal{D}_{h}'\left(\mathbb{R}^{n}\right)\rightarrow
C^{\infty}\left(\mathbb{R}^{n}\right).
\end{equation*}
We shall also use the notation $T\equiv T'.$

Lastly, we define global semi-classical Fourier integral operators.
\begin{Def}\label{dfio}
Let $M$ be a smooth $k$-dimensional manifold and let $\Lambda\subset
T^{*}M$ be a smooth closed
Lagrangian submanifold with respect to the canonical symplectic 
structure on $T^{*}M.$
Let $r\in\mathbb{R}.$
Then the space $I^{r}_{h}\left(M, \Lambda\right)$ of semi-classical
Fourier integral
distributions of order $r$ associated to $\Lambda$ is defined as the set
of all $u\in\mathcal{D}'_{h}\left(M\right)$ 
such
that
\begin{equation}\label{defgfio}
\left(\prod_{j=0}^{N}
A_{j}\right)\left(u\right)=\mathcal{O}_{L^{2}\left(M\right)}\left(h^{N-r-\frac{k}{4}}\right),
h\to 0,
\end{equation}
for all $N\in\mathbb{N}_{0}$ and for all $A_{j}\in \Psi_{h}^{0}\left(1,
X\right),$ $j=0, \dots, N-1,$ with
compactly
supported symbols and principal symbols vanishing on $\Lambda$, and any $ 
A_N \in 
\Psi_h^{0} ( 1 , X ) $ with a compactly supported symbol.

A continuous linear operator
$C_{c}^{\infty}\left(M_1\right)\rightarrow\mathcal{D}_{h}'\left(M_2\right),$ 
where $M_1, M_2$ are smooth manifolds, 
whose Schwartz kernel is an element of 
$I_{h}^{r}(M_1\times M_2, \Lambda)$ for some 
Lagrangian submanifold $\Lambda\subset T^{*}M_1\times T^{*}M_2$ and some $r\in\mathbb{R}$ 
will be called a global semi-classical Fourier integral
operator of order $r$ associated to $\Lambda.$
We denote the space of these operators by 
$\mathcal{I}_{h}^{r}(M_1\times M_2, \Lambda).$
\end{Def}

\subsection{Representation of the Scattering Amplitude}\label{sreprsa}
In this section we define the semi-classical scattering amplitude and
derive the representation of the scattering amplitude, which
will be used in the proof of Theorem 1.
The derivation is similar to the one presented in \cite{S}.

We recall that for any $\theta\in \mathbb{S}^{n-1}$ and
$\lambda>0$ there exists a unique, up to a compactly supported function, solution $\psi$ to the problem
$\left(P\left(h\right)-\lambda\right)\psi(\cdot, \theta; \lambda, h)=0$,
$\psi(\cdot, \theta; \lambda,
h)\in\mathcal{D}_{loc}\left(P\left(h\right)\right)$ such that
\begin{equation*}
\psi\left(x, \theta; \lambda,
h\right)=\left(1-\chi\left(x\right)\right)e^{\frac{i\sqrt{2\lambda}\theta\cdot
x}{h}}+\psi_{sc}(x, \theta; \lambda, h),
\end{equation*}
where $\psi_{sc}$ satisfies the Sommerfeld outgoing condition at 
infinity:
$$\left(\partial/\partial r
-i\sqrt{2\lambda}/h\right)\psi_{sc}=\mathcal{O}\left(r^{-\left(n+1\right)/2}\right), 
\text{ as } r=\|x\|\to
\infty$$ and $\chi\in C_{c}^{\infty}\left(\mathbb{R}^{n}\right)$ is 
equal
to 1 on
$B\left(0, R_{0}\right).$
Then
\begin{equation*}
\psi\left(x, \theta; \lambda,
h\right)=e^{\frac{i\sqrt{2\lambda}\theta\cdot
x}{h}}+\frac{e^{i\sqrt{2\lambda}r/h}}{r^{\left(n-1\right)/2}}A\left(\frac{x}{r},
\theta; \lambda,
h\right)+\mathcal{O}\left(\frac{1}{r^{\left(n+1\right)/2}}\right),
\text{as } r=\|x\|\to\infty.
\end{equation*}
The function $A$ is called the scattering amplitude.

Let, now, $\chi_{1}\in C_{c}^{\infty}\left(\mathbb{R}^{n}\right)$ be equal
to 1
on $B\left(0,
R_{0}\right)$ and let
$\chi_{2}\in C_{c}^{\infty}\left(\mathbb{R}^{n}\right)$ be equal to 1 the
support of
$\chi_{1}$.
Let
\begin{displaymath}
\psi=\left(1-\chi_{1}\right)e^{\frac{i\sqrt{2\lambda}\left\langle \theta,
\cdot\right\rangle}{h}}+\psi_{sc}
\end{displaymath}
be such that $$\left(P\left(h\right)-\lambda\right)\psi=0.$$
Then 
\begin{equation}\label{eq}
\begin{aligned}
\left(P\left(h\right)-\lambda\right)\psi_{sc}
&=-\left(P\left(h\right)-\lambda\right)(1-\chi_{1})e^{\frac{i\sqrt{2\lambda}\left\langle 
\theta, \cdot \right\rangle}{h}}\\
&=-\left(-\frac{1}{2}h^{2}\Delta-\lambda\right)\left(1-\chi_{1}\right)e^{\frac{i\sqrt{2\lambda}\left\langle 
\theta, 
\cdot\right\rangle}{h}}=-\frac{1}{2}\left[h^{2}\Delta,\chi_{1}\right]
e^{\frac{i\sqrt{2\lambda}\left\langle \theta, \cdot \right\rangle}{h}}
\end{aligned}
\end{equation}
and therefore 
\begin{equation}\label{1stpsi}
\psi_{sc}=-\frac{1}{2}R\left(\lambda, 
h\right)\left[h^{2}\Delta,\chi_{1}\right]e^{\frac{i\sqrt{2\lambda}\left\langle 
\theta, \cdot \right\rangle}{h}}.
\end{equation}
Then 
\begin{equation}\label{2ndpsi1}
\left(-\frac{1}{2}h^{2}\Delta-\lambda\right)\left(1-\chi_{2}\right)\psi_{sc}=\frac{1}{2}\left[h^{2}\Delta,\chi_{2}\right]\psi_{sc}
+\left(1-\chi_{2}\right)\left(P\left(h\right)-\lambda\right)\psi_{sc}
\end{equation}
Substituting (\ref{eq}) into (\ref{2ndpsi1}), we 
obtain
\begin{equation}\label{2ndpsi2} 
\begin{aligned}
\left(-\frac{1}{2}h^{2}\Delta-\lambda\right)\left(1-\chi_{2}\right)\psi_{sc}
&=\frac{1}{2}\left[h^{2}\Delta,\chi_{2}\right]\psi_{sc}
-\frac{1}{2}\left(1-\chi_{2}\right)\left[h^{2}\Delta,\chi_{1}\right]e^{\frac{i\sqrt{2\lambda}\left\langle\theta, 
\cdot\right\rangle}{h}}\\
 & =\frac{1}{2}\left[h^{2}\Delta,\chi_{2}\right]\psi_{sc}
\end{aligned}
\end{equation}
since $\supp \nabla \chi_{1}\cap \supp \left(1-\chi_{2}\right)=\emptyset.$
Substituting (\ref{1stpsi}) into (\ref{2ndpsi2}), we 
obtain
\begin{equation*}
\left(1-\chi_{2}\right)\psi_{sc} 
=-\frac{1}{4}R_{0}\left(z,h\right)\left[h^{2}\Delta, 
\chi_{2}\right]R\left(z,h\right)\left[h^{2}\Delta, 
\chi_{1}\right]e^{\frac{i\sqrt{2\lambda}\left\langle \theta, 
\cdot\right\rangle}{h}}.
\end{equation*}
Therefore, by Proposition 1.1 in \cite{Mel}, we obtain that the scattering 
amplitude $A$ is 
\begin{equation*}
A\left(\omega, \theta; \lambda, h\right)=c\left(n, \lambda; h\right) \int 
e^{-\frac{i\sqrt{2\lambda}\left\langle 
\omega, x \right\rangle}{h}}\left[h^{2}\Delta, \chi_{2}\right]R\left(\lambda, 
h\right)\left[h^{2}\Delta, 
\chi_{1}\right]e^{\frac{i\sqrt{2\lambda}\left\langle \theta, 
\cdot\right\rangle}{h}}dx, 
\end{equation*}
where $c\left(n, \lambda, 
h\right)=\frac{e^{-i\pi\frac{n-3}{4}}\lambda^{\frac{n-3}{4}}}
{2^{\frac{n+13}{4}}\left(\pi h\right)^{\frac{n+1}{2}}}$
and this is the representation of the scattering amplitude, which we will 
use in the proof of the main theorem.
The independence of this representation of the choice of the functions 
$\chi_{j}$, $j=1, 2$ with the above properties is proved in \cite{PZ}.
Proposition 2.1 in \cite{PZ} further shows that the scattering amplitude 
with the constant $c(n, \lambda, h)$ replaced by $\tilde{c}(n, \lambda, 
h)=\frac{i\pi\lambda^{\frac{n-2}{2}}}{(2\pi)^{n}h^{n}}$ is also
the kernel of $S(\lambda, h)-I,$ where $S(\lambda, h)$ is the scattering 
matrix at energy $\lambda.$ 

From this representation of the scattering amplitude it is clear that it 
can be extended meromorphically everywhere where the resolvent can be 
extended meromorphically and that the poles of the scattering amplitude are among the resonances.

We shall now prove two lemmas which give further information about the 
structure of the cut-off resolvent and the scattering amplitude.

\begin{Lem}\label{rsym}
For $\lambda>0$ not a resonance $K_{R(\lambda, h)}\left(x, y\right)=K_{R(\lambda,
h)}(y, x)$ for $x, y\in
\mathbb{R}^{n}\backslash B\left(0, R_{0}\right).$
\end{Lem}

\medskip
\noindent
{\it Proof:}\;
For $u, v\in L^{2}\left(\mathbb{R}^{n}\backslash B\left(0, 
R_{0}\right)\right)$ let
$\left\langle u, v\right\rangle=\int u v.$
Let $u$ and $v$ be further chosen with compact support and
let $z\in\mathbb{C}$ be such that $\Im z>0.$
We then have
\begin{equation}\label{symR}
\begin{aligned}
& \left\langle \left(R\left(z, h\right)u\right)|_{\mathbb{R}^{n}\backslash
B\left(0, R_{0}\right)}, v\right\rangle
=\left\langle \left(R\left(z, h\right)u\right)|_{\mathbb{R}^{n}\backslash
B\left(0, R_{0}\right)}, \left(P\left(h\right)-z\right)R\left(z,
h\right)v\right\rangle\\
& \quad =\left\langle \left(R\left(z,
h\right)u\right)|_{\mathbb{R}^{n}\backslash B\left(0, R_{0}\right)},
\left(-\frac{1}{2}h^{2}\Delta-z\right)\left(\left(R\left(z,
h\right)v\right)|_{\mathbb{R}^{n}\backslash B\left(0,
R_{0}\right)}\right)\right\rangle\\
& \quad =\left\langle
\left(-\frac{1}{2}h^{2}\Delta-z\right)\left(\left(R\left(z,
h\right)u\right)|_{\mathbb{R}^{n}\backslash B\left(0,
R_{0}\right)}\right), \left(R\left(z,
h\right)v\right)|_{\mathbb{R}^{n}\backslash B\left(0,
R_{0}\right)}\right\rangle\\
& \quad =\left\langle \left(\left(P\left(h\right)-z\right)\left(R\left(z,
h\right)u\right)\right)|_{\mathbb{R}^{n}\backslash B\left(0, R_{0}\right)},
\left(R\left(z, h\right)v\right)|_{\mathbb{R}^{n}\backslash B\left(0,
R_{0}\right)}\right\rangle\\
& \quad =\left\langle u, \left(R\left(z,
h\right)v\right)|_{\mathbb{R}^{n}\backslash B\left(0,
R_{0}\right)}\right\rangle.
\end{aligned}
\end{equation}

Let, now, $\lambda\in\mathbb{R}\backslash\left\{0\right\}$ and let
$\left(z_{k}\right)_{k\in\mathbb{N}}\subset\mathbb{C}$ satisfy $\Im z_{n}>0,$
$z_{k}\to\lambda,$ $k\to\infty.$
Then, from (\ref{symR}) we have that for every $k$
\begin{equation}\label{zk}
\left\langle \left(R\left(z_{k},
h\right)u\right)|_{\mathbb{R}^{n}\backslash B\left(0, R_{0}\right)},
v\right\rangle
=\left\langle u, \left(R\left(z_{k},
h\right)v\right)|_{\mathbb{R}^{n}\backslash B\left(0,
R_{0}\right)}\right\rangle.
\end{equation}

Letting $k\to\infty$ in (\ref{zk}) and using the fact that
$R\left(\cdot,
h\right):\mathcal{H}_{comp}\rightarrow\mathcal{H}_{loc}$ is analytic in
the
upper-half plane and up to the real axis, we obtain
\begin{equation*}
\left\langle \left(R\left(\lambda,
h\right)u\right)|_{\mathbb{R}^{n}\backslash B\left(0, R_{0}\right)},
v\right\rangle
=\left\langle u, \left(R\left(\lambda,
h\right)v\right)|_{\mathbb{R}^{n}\backslash B\left(0, R_{0}\right)}\right\rangle,
\end{equation*}
which completes the proof of the lemma.$\hfill\Box$

\begin{Lem}
The operator $\psi_1 R\left(z, h\right)\psi_2$ has a smooth 
Schwartz kernel 
$K_{\psi_1 R\left(z, h\right)\psi_2}$ when $z$ is not a resonance and 
$\psi_j\in 
C_{c}^{\infty}(\mathbb{R}^{n}\backslash B\left(0, R_{0}\right)), j=1, 2,$ have 
disjoint 
supports.
\end{Lem}

\medskip
\noindent
{\it Proof:}\;
Let $\varphi\in C^{\infty}\left(\mathbb{R}^{n}\right)$ and 
consider
\begin{equation*}
\begin{aligned}
& \left(-\frac{1}{2}h^{2}\Delta-z\right) \psi_1 R\left(z, 
h\right)\psi_2\varphi\\
& =\left[-\frac{1}{2}h^{2}\Delta, \psi_1\right]R\left(z, 
h\right)\psi_2\varphi+\psi_1\left(-\frac{1}{2}h^{2}\Delta-z\right)\left(\left(R\left(z, 
h\right)\psi_2\varphi\right)|_{\mathbb{R}^{n}\backslash B\left(0, 
R_{0}\right)}\right)\\
 & =\left[-\frac{1}{2}h^{2}\Delta, \psi_1\right]R\left(z, 
h\right)\psi_2\varphi+\psi_1\left(\left(P\left(h\right)-z\right)R\left(z, 
h\right)\psi_2\varphi\right)|_{\mathbb{R}^{n}\backslash B\left(0, 
R_{0}\right)}\\
 & =\left[-\frac{1}{2}h^{2}\Delta, \psi_1\right]R\left(z, 
h\right)\psi_2\varphi+\psi_1\psi_2\varphi\\
 & =\left[-\frac{1}{2}h^{2}\Delta, \psi_1\right]R\left(z, 
h\right)\psi_2\varphi\in 
H^{-1}(\mathbb{R}^{n}).
\end{aligned}
\end{equation*}
Therefore,
\begin{equation*}
\psi_1 R\left(z, h\right)\psi_2\varphi\in H^{1}(\mathbb{R}^{n}).
\end{equation*}
Similarly, for every $k\in\mathbb{N}$ we have that
\begin{equation*}
\begin{aligned}
& \left(-\frac{1}{2}h^2 \Delta-z\right)^k\psi_1 R\left(z, 
h\right)\psi_2\varphi\\
& \quad = \underbrace{\left[-\frac{1}{2}h^2\Delta,\left[\dots,\left[-\frac{1}{2}h^2\Delta, 
\psi_1\right]\dots\right]\right]}_k 
R\left(z, 
h\right)\psi_2\varphi\in H^{-k}(\mathbb{R}^{n}) 
\end{aligned}
\end{equation*}
and therefore
\begin{equation*}
\psi_1 R\left(z, h\right)\psi_2\varphi\in H^{k}(\mathbb{R}^{n})\text{ for 
every } k.
\end{equation*}
Thus
\begin{equation}\label{smooth}
\psi_1 R\left(z, h\right)\psi_2\varphi\in 
C^{\infty}\left(\mathbb{R}^n\right).
\end{equation}
This, together with Lemma \ref{rsym}, implies that the kernel of 
$\psi_{1}R\left(\lambda, h\right)\psi_{2}$ is smooth.$\hfill\Box$

As in \cite{PZ} we, now, introduce the operators
\begin{equation*}
\left[\mathbb{E}_{\pm}\left(\lambda, 
h\right)f\right]\left(\omega\right)=\int 
e^{\pm\frac{i\sqrt{2\lambda}\omega\cdot 
x}{h}}f\left(x\right) dx = 
\hat{f}\left(\mp \frac{\sqrt{2\lambda}\omega}{h}\right), \omega\in 
S^{n-1},
\end{equation*}
where $f$ has compact support.
Then we can express $A$ as
\begin{equation}\label{amplk}
\begin{aligned}
A\left(\lambda, h\right) & =c\left(n, \lambda, h\right)\left(\mathbb{E}_{-}\left(\lambda, 
h\right)\otimes \mathbb{E}_{+}\left(\lambda, 
h\right)\right)\left(\left[h^{2}\Delta, \chi_{2}\right]\otimes 
\left[h^{2}\Delta, 
\chi_{1}\right]^{t}\right)K_{\tilde{\chi}_2 R(z, h)\tilde{\chi}_1},
\end{aligned}
\end{equation}
where $\tilde{\chi}_{j}\in C_{c}^{\infty}(\mathbb{R}^{n}\backslash B(0, 
R_{0}))$ are such that $\tilde{\chi}_{j}=1$ on $\supp\nabla\chi_{j},$ $j=1, 
2$ and 
$\supp\tilde{\chi}_{1}\cap\supp\tilde{\chi}_{2}=\emptyset.$$\hfill\Box$

\section{Scattering Geometry}\label{cgeom}
Here we collect the geometric results which we will use in the proof of 
the main theorem and in the applications. 

\subsection{Scattering Relation}\label{scr}
Here we define the canonical relation, which we will prove to be quantized
by the scattering amplitude in the sense of semi-classical Fourier
integral operators.
We shall work in the following setting.
Let $X$ be a smooth manifold of dimension $n>1$ such that
$X$ coincides with $\mathbb{R}^{n}$ outside of $B(0, R_{0})$ for some $R_{0}>0.$
Let $g$ be a Riemannian metric on $X$ which satisfies the condition
\begin{equation*}
g_{ij}\left(x\right)=\delta_{ij} \text{ for } \|x\|>R_{0}.
\end{equation*}
Let $V\in C_{c}^{\infty}\left(X\backslash B(0, R_{0})^{c}; \mathbb{R}\right).$
Let $P\left(h\right)=\frac{1}{2}h^{2}\Delta_{g}+V,$ $0<h\leq 1,$ with $p(x, 
\xi)=\frac{1}{2}\|\xi\|_{g}+V(x)$ denoting its semi-classical principal symbol. 
We assume that for some $\lambda>0$ the operator 
$P(h)-\lambda$ is of princiapl type.
This implies that $\Sigma_{\lambda}=p^{-1}\left(\lambda\right)$ is a
smooth $2n-1$
dimensional manifold.
Let $(x_{0}, \xi_{0})\in\Sigma_{\lambda}$ with
$x_{0}\notin B(0, R_{0})^{c}$ be such that there exists $T>0$ 
satisfying $x(t; x_{0}, \xi_{0})\in
B(0, R_{0})^{c}$ for all
$|t|>T.$
The curve $\gamma(x_{0}, \xi_{0})$ is called a non-trapped phase trajectory.
Then it follows that there exists an open neighborhood $W\subset T^{*}X$ 
of $(x_{0},
\xi_{0})$ such that for all $(x', \xi')\in W\cap \Sigma_{\lambda}$ we have 
that $x(t; x',
\xi')\in B(0, R_{0})^{c}$ for all $|t|>T.$

Let, now, $i: T^{*}\mathbb{S}^{n-1}\hookrightarrow T^{*}\mathbb{R}^{n}$ 
denote the
inclusion map.
Let $\psi: T^{*}\mathbb{R}^{n}\ni \left(x, \xi\right)\mapsto \left(\xi,
x\right)\in
T^{*}\left(\mathbb{R}^{n}\right)^{*}$.
Then $L=\left(\psi\circ i\right) \left(T^{*}\mathbb{S}^{n-1}\right)$ is a
smooth submanifold of
$T^{*}\mathbb{R}^{n}$ and therefore
\begin{equation*}
L_{1}\left(\lambda\right)=\left\{\left(x, \xi\right)\bigg|
\left(x+\frac{\left(1+R_{0}\right)}{\sqrt{2\lambda}}\xi,
\frac{\xi}{\sqrt{2\lambda}}\right)\in
L\right\}\subset \Sigma_{\lambda}
\end{equation*}
and
\begin{equation*}
L_{2}\left(\lambda\right)=\left\{\left(x, \xi\right)\bigg|
\left(x+\frac{\left(1-R_{0}\right)}{\sqrt{2\lambda}}\xi,
\frac{\xi}{\sqrt{2\lambda}} \right)
\in L\right\}\subset \Sigma_{\lambda}
\end{equation*}
are hypersurfaces transverse to $H_{p}$ near $\gamma(\cdot; x_{0}, 
\xi_{0}).$

Since $\gamma(\cdot; x_{0}, \xi_{0})$ is a non-trapped phase trajectory 
and since
the support of the perturbation
is compact, the Hamiltonian flow of $p$ on $\Sigma_{\lambda}$ near 
$\gamma(\cdot; x_{0},
\xi_{0})$ outside the
support of the
perturbation is
the free Hamiltonian flow defined by the equations
\begin{equation*}
\left\{\begin{array}{l}
\dot{\xi}=0\\
\dot{x}=\xi
\end{array}\right.
\end{equation*}
with solution
\begin{equation*}
\left\{\begin{array}{l}
\xi=\sqrt{2\lambda}\xi_{\infty}=\text{const.}\\
x=\sqrt{2\lambda}\xi_{\infty}t+x_{\infty}
\end{array}\right.
\end{equation*}
where by conservation of energy, we have that
$\xi_{\infty}\in\mathbb{S}^{n-1}.$

\begin{figure}[t]
\begin{center}
\input{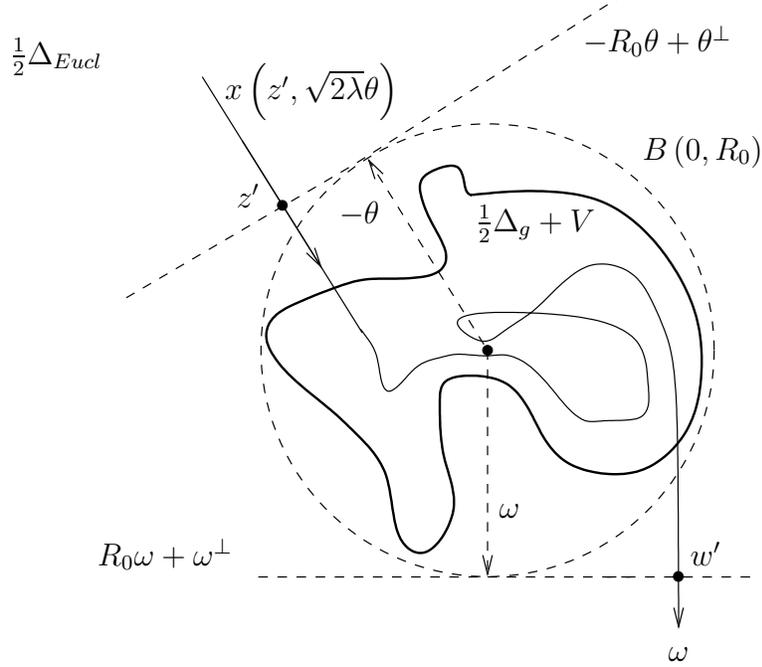}
\end{center}
\caption{The scattering relation consists of the points
$\left(\sqrt{2\lambda}\theta, z';
\sqrt{2\lambda}\omega, w'\right)$ as in this figure.}
\label{fig:sr}
\end{figure}

Denoting $p_1=\gamma(\cdot; x_{0}, \xi_{0})\cap L_{1}(\lambda)$ we then 
have that
for
every $\left(x, \xi\right)\in L_{1}\left(\lambda\right)$ in a sufficiently 
small
neighborhood of $p_{1}$ there exists a
unique $T\left(x, \xi\right)>0$ such that 
\[\exp \left(T\left(x, \xi\right)H_{p}\right)\left(x, \xi\right)\in
L_{2}\left(\lambda\right).\]
Therefore, since $L_{1}\left(\lambda\right)$ and
$L_{2}\left(\lambda\right)$ are hypersurfaces in $\Sigma_{\lambda}$ transverse to $H_p$ near 
$\gamma,$ we have that there exists an open set $U\subset L_{1}(\lambda),$ $p_{1}\in U,$ such 
that
\begin{equation*}
\begin{split}
SR_U\left(\lambda\right)= & \bigg\{\left(x, \xi; y,
\eta\right)\bigg|\left(z,
\zeta\right)=\left(\frac{x}{\sqrt{2\lambda}}-\left(R_{0}+1\right)\xi,
\sqrt{2\lambda}\xi\right)\in U,\\
& \left(\frac{y}{\sqrt{2\lambda}}+\left(R_{0}-1\right)\eta,
\sqrt{2\lambda}\eta\right)=
\exp \left(T\left(z,
\zeta\right)H_{p}\right)\left(z, \zeta\right)\bigg\}'
\end{split}
\end{equation*}
is a Lagrangian submanifold of $\left(T^{*}\mathbb{S}^{n-1}\times
T^{*}\mathbb{S}^{n-1},
\pi_{1}^{*}\sigma+\pi_{2}^{*}\sigma\right),$ where $\pi_{j}:
T^{*}\mathbb{S}^{n-1}\times
T^{*}\mathbb{S}^{n-1}\to T^{*}\mathbb{S}^{n-1},$ $j=1, 2,$ is the 
projection onto the $j$-th
factor.
We shall call $SR_U\left(\lambda\right)$ the scattering relation at energy
$\lambda$ (see Figure \ref{fig:sr}). 

We now show how, under a certain geometric assumption, we can find a phase 
function which
parameterizes the scattering relation near a non-trapped trajectory.
To state the assumption, let us first introduce some notation.
For $\theta\in\mathbb{S}^{n-1}$ and $z\in\theta^{\perp}-R_0\theta,$ if 
$\gamma\left(\cdot; z, \sqrt{2\lambda}\theta\right)$ is a  non-trapped
trajectory, then, as we saw above, there exist
\[x_{\infty}\left(\theta, z\right), \xi_{\infty}\left(\theta, z\right)\in\mathbb{R}^{n},\;\left\|\xi_{\infty}\left(\theta, z\right)\right\|=1,\] such that
\begin{equation*}
\gamma\left(t; z, \sqrt{2\lambda}\theta\right)=\left(x_{\infty}\left(\theta,
z\right)+t\sqrt{2\lambda}\xi_{\infty}\left(\theta, z\right),
\sqrt{2\lambda}\xi_{\infty}\left(\theta, z\right)\right),
t> > 0.
\end{equation*}
We shall call such a non-trapped phase trajectory with initial direction
$\theta$ and final
direction
$\omega=\xi_{\infty}\left(\theta, z\right),$ a $(\theta, 
\omega)$-trajectory.
We,  now, make the following

\begin{Def}
If $\theta_{0}, \omega_{0}\in\mathbb{S}^{n-1}$ are such that
\begin{equation*}
\xi_{\infty}(\theta_{0}, z')=\omega_{0}
\end{equation*}
implies that the map
\begin{equation}\label{regular}
\theta_{0}^{\perp}-R_0\omega_0\ni z\mapsto \xi_{\infty}\left(\theta_{0},
z\right)\in\mathbb{S}^{n-1}
\text{ is non-degenerate at } z',
\end{equation}
then we shall say that $\omega_{0}$ is regular for $\theta_{0}.$
\end{Def}
We remark that this definition is a rephrasing of condition \eqref{angulard}.

We then have the following

\begin{Lem}\label{regL}
If $\omega_0\in\mathbb{S}^{n-1}$ is regular for 
$\theta_0\in\mathbb{S}^{n-1},$ then
\begin{enumerate}[(a)]
\item $\theta_{0}\ne\omega_{0}.$
\item There exist $O_j\subset\mathbb{S}^{n-1},$ $j=1, 2,$ open, 
$\theta_0\in 
O_1,$ $\omega_0\in
O_2,$  and a number $L\in\mathbb{N}$ such that for every
$(\theta, \omega)\in O_1\times O_2,$ there exist at least 
$L$ $(\theta, \omega)-$ trajectories.
\end{enumerate}
\end{Lem}

\medskip
\noindent
{\it Proof:}\;
We shall work in a local trivialization of $T^{*}\mathbb{S}^{n-1}_{R_{0}}$ near 
$T^{*}_{R_{0}\theta_{0}}\mathbb{S}^{n-1}_{R_{0}}$ and 
$T^{*}_{R_{0}\omega_{0}}\mathbb{S}^{n-1}_{R_{0}},$ where $\mathbb{S}^{n-1}_{R_{0}}=\{x\in\mathbb{R}^{n}: \|x\|=R_0\}.$
We first prove that $\theta_{0}\ne\omega_{0}.$
Assume that $\theta_0=\omega_0.$
Then for every $z\in\mathbb{R}^{n-1}$ with
$\|z\|>R_{0}$ we have that $\xi_{\infty}(\theta_{0}, z)=\theta_{0}.$
Therefore
$\det \left(\frac{\partial \xi_{\infty}(\theta_{0}, \cdot)}{\partial
z}(z)\right)=0$ for $z\in\theta_{0}^{\perp},$ $\|z\|>R_{0},$ which is a
contradiction with
the regularity assumption.
Thus it follows that $\theta_{0}\ne\omega_{0},$ which establishes (a).

Let, now, $z'\in\left(\xi_{\infty}(\theta_0, \cdot)\right)^{-1}(\omega_{0}).$
Then, by the Inverse Function Theorem, there exist open sets 
$O'_{2}\subset\mathbb{S}^{n-1},$ $\omega_{0}\in O'_{2}$ and 
$\mathcal{Z}'\subset\mathbb{R}^{n-1},$ $z'\in\mathcal{Z}'$ such that 
$\xi_{\infty}|_{\mathcal{Z}'}\left(\theta_{0}, \cdot\right)$ is a diffeomorphism onto 
$O'_{2}.$
Therefore, the set $\left(\xi_{\infty}(\theta_0, \cdot)\right)^{-1}(\omega_{0})$ is discrete.
From the first part of this proof, it follows that 
$\left(\xi_{\infty}(\theta_0, \cdot)\right)^{-1}(\omega_{0})$ is also bounded.
Therefore, it is finite and we shall denote its elements by $\{z_1, \dots, z_L\},$ $L\in\mathbb{N}.$ 

By the Implicit Function Theorem and the regularity assumption, we now have
that there exist open sets $O_1, O_2\subset\mathbb{S}^{n-1}$ with
$\theta_{0}\in O_1$ and $\omega_{0}\in O_2$ and functions $z_{l}\in
C^{\infty}(O_1\times O_2; \mathbb{R}^{n-1}),$ $l=1, 
\dots, L,$ such that
$z_{l}(\theta_0, \omega_0)=z_l$ and $\xi_{\infty}(\theta, z_{l}(\theta,
\omega))=\omega,$ $(\theta, \omega)\in O_1\times O_2,$ which completes the proof of 
(b).$\hfill\Box$

Let, now, $w_{l}(\theta, \omega)=\gamma\left(\cdot; z_l\left(\theta, \omega\right), \sqrt{2\lambda}\theta\right)\cap \left(R_{0}\omega+\omega^{\perp}\right).$
Then, as above we have that, 
\begin{equation}
\begin{aligned}
SR^{l}_{O_{1}\times O_{2}}(\lambda)=\Big\{\Big(\sqrt{2\lambda}\theta, 
\sqrt{2\lambda}\omega, z_{l}(\theta, \omega)+(R_0-1)\theta, w_{l}(\theta, 
\omega)&-(R_0-1)\omega\Big)\Big|\\
& (\theta, \omega)\in O_1\times O_2\Big\}'
\end{aligned}
\end{equation}
is a Lagrangian submanifold of $T^{*}(\mathbb{S}^{n-1}\times\mathbb{S}^{n-1}).$Furthermore,  
$\pi|_{SR_{O_{1}\times O_{2}}^{l}(\lambda)}$ is a surjection, where 
$\pi:T^{*}(\mathbb{S}^{n-1}\times\mathbb{S}^{n-1})\to 
\mathbb{S}^{n-1}\times\mathbb{S}^{n-1}$ is the canonical projection.
Therefore, after decreasing, if necessary, $O_1\times O_2$ around $\left(\theta_0, \omega_0\right)$ we have that there exists a function $\mathcal{W}_{l}\in 
C^{\infty}(O_1\times O_2)$ such that 
\[SR_{O_{1}\times O_{2}}^{l}(\lambda)=\left\{(m, d\mathcal{W}_{l}(m)): 
m\in O_1\times O_2\right\}.\] 
Since we also have that $SR_{O_{1}\times O_{2}}^{l}(\lambda)$ is a canonical relation and therefore, after possibly decreasing 
$O_1\times O_2$ further near $(\theta_0, \omega_0),$ we can assume that 
\begin{equation}\label{detW}
\det\left(\partial_{\theta} 
\partial_{\omega}\mathcal{W}_{l}(\theta, 
\omega)\right)\ne 0,\,
\left(\theta, \omega\right)\in O_1\times O_2.
\end{equation}
We let 
\begin{equation}\label{lagrl}
SR_{l}(\lambda)=SR^{l}_{O_{1}\times O_{2}}(\lambda)
\end{equation}
with $O_{1}\times O_{2}$ as in \eqref{detW}.
We observe that \eqref{detW} implies that the map $R_{0}\omega+\omega^{\perp}\ni w\mapsto 
\xi_{\infty}\left(-\omega, w\right)\in \mathbb{S}^{n-1}$ is non-degenerate in 
a neighborhood of $w_{l}(\theta, \omega),$ $\left(\theta, \omega\right)\in O_1\times O_2.$

For the $(\theta, \omega)$ trajectory defined by $z_{l}(\theta, \omega)$ as in the proof of 
the Lemma, we shall use the subscript $l$ to distinguish it from all other $(\theta, 
\omega)$ trajectories.

The same proof as in \cite[Lemma 3.2]{RT} together with \eqref{detW} show that there 
exist
$0<S_{0}<S_{1}$ and $T_0>>0$ and open sets $U^{l}_{\theta, 
\omega}\subset\mathbb{R}^{n-1},$ $z_{l}(\theta, \omega)\in U^{l}_{\theta, 
\omega},$ 
$(\theta, \omega)\in O_1\times O_2,$ $l=1, \dots, L$ such that
\begin{equation}\label{candefaction}
\det\left(\frac{\partial x_{l}\left(t; \cdot, 
\sqrt{2\lambda}\theta\right)}{\partial
y}\left(y\right)\right)\ne 0
\end{equation}
for $\theta\in O_1,$ $y\in\left\{z-s\theta: z\in U^{l}_{\theta, \omega},
s\in\left[S_0, S_1\right]\right\},$ $\omega\in O_2,$ $t>T_0.$

For $(\theta, \omega)\in O_1\times O_2$ we now define the (modified) action along the 
segment of the $l-$th $\left(\theta, \omega\right)$-trajectory, 
$\gamma\left(\cdot; z_l(\theta, \omega), \sqrt{2\lambda}\theta\right),$ 
between the points $y_{l}\left(s; \theta, \omega\right)=z_{l}\left(\theta,
\omega\right)-\sqrt{2\lambda}s\theta\in
W_{\theta},$ $s\in[S_1, S_0]$ and $x_{l}(t, s, \theta, \omega)=x_{l}\left(t; y_l(s;
\theta, \omega), \sqrt{2\lambda}\theta\right),$ $t>T_0.$
We choose a fixed $t>T_0$ and we set
\begin{equation}\label{maction}
S_{l}\left(\theta, \omega\right)=\left\langle y_{l}\left(s; \theta, \omega\right),
\sqrt{2\lambda}\theta\right\rangle+\int_{0}^{t} L\left(x_{l}, 
\dot{x}_{l}\right)dt-\left\langle x_l\left(t, s, \theta, \omega\right),
\sqrt{2\lambda}\omega\right\rangle-\lambda t,
\end{equation}
where $L(x, \dot{x})=\frac{1}{2}\left\|\dot{x}\right\|^{2}_{g}-V(x)$ is 
the Lagrangian, and the integral is taken over the $l-$th bicharacteristic curve
connecting $y_{l}\left(s; \theta, \omega\right)$ and $x_{l}\left(t, s, \theta, 
\omega\right).$
We observe that, since the support of the perturbation is compact,
$S_l(\theta, \omega)$ is independent of $s$ for
$s\in[S_1, S_0].$

\begin{Lem}\label{actionsr}
Let $\omega_{0}\in\mathbb{S}^{n-1}$ be regular for 
$\theta_{0}\in\mathbb{S}^{n-1}.$

Then $SR_{l}(\lambda)=\Lambda_{S_{l}},$ where $\Lambda_{S_{l}}=\left\{\left(\theta, \omega, d_{\theta}S_{l}, d_{\omega}S_{l}\right): (\theta, \omega)\in\mathbb{S}^{n-1}\times\mathbb{S}^{n-1}\right\},$ $l=1, \dots, L.$
\end{Lem}

\medskip
\noindent
{\it Proof:}\;

We consider
\begin{equation}\label{somega}
\begin{aligned}
d_{\omega}S_{l}(\theta, \omega) & =d_{\omega}\left(\left\langle 
y_l\left(s; \theta, \omega\right),
\sqrt{2\lambda}\theta\right\rangle+\int_{0}^{t} L\left(x_l,
\dot{x}_l\right)dt\right)-d_{\omega}\left(\left\langle x_l\left(t, s, \theta, 
\cdot\right),
\sqrt{2\lambda}\cdot\right\rangle \right)(\omega)\\
& =\left\langle\sqrt{2\lambda}\omega, d_{\omega}x_l(t, s, \theta, 
\cdot)(\omega)\right\rangle
-\left\langle\sqrt{2\lambda}\omega, d_{\omega}x_l(t, s, \theta, 
\cdot)(\omega)\right\rangle\\
& \quad\quad -d_{\omega}\left\langle x_l(t, s, \theta, \omega), 
\sqrt{2\lambda}\cdot\right\rangle(\omega)\\
& = -d_{\omega}\left\langle x_l\left(\omega, \theta, s, t\right),
\sqrt{2\lambda}\cdot\right\rangle(\omega),
\end{aligned}
\end{equation}
where \eqref{candefaction} has allowed us to use \cite[Theorem 46.C]{A} 
to obtain the second equality.

To compute $d_{\theta}S_l$ we first reparameterize the phase
trajectories in the reverse direction, which is equivalent to
considering the reverse of the initial and final directions.
We further re-write $S_l\left(\theta, \omega\right)$ in the 
following way
\begin{equation*}
S_l\left(\theta, \omega\right)=-\left\langle x_l\left(s; \theta,
\omega\right), \sqrt{2\lambda}\omega\right\rangle +\int_{0}^{t} L\left(x_l,
\dot{x}_l\right)dt+\left\langle y_l\left(t, s, \omega, \theta\right),
\sqrt{2\lambda}\theta\right\rangle-\lambda t,
\end{equation*}
where $x_l\left(s; \theta, \omega\right)=w_{l}\left(\theta, 
\omega\right)+\sqrt{2\lambda}\omega s,$
$y_l\left(t, s, \theta, \omega\right)=x_l\left(t;
x_l\left(s; \theta, \omega\right),-\sqrt{2\lambda}\omega\right),$
$s\in[S_0, S_1],$ and the integral is taken over the bicharacteristic curve connecting 
$x_l\left(s; \theta, \omega\right)$ and $y_l\left(t, s, \theta, \omega\right).$
We observe that this bicharacteristic curve is uniquely defined by \eqref{detW} and \eqref{candefaction}.

Equations \eqref{detW} and \eqref{candefaction} further allow us to proceed as in (\ref{somega}) and we obtain
\begin{equation}\label{stheta}
\begin{aligned}
d_{\theta}S_l\left(\omega,
\theta\right) & = d_{\theta}\left(-\left\langle x_l\left(s; \theta,
\omega\right),
\sqrt{2\lambda}\omega\right\rangle +\int_{0}^{t} L\left(x_l,
\dot{x}_l\right)dt\right)+d_{\theta}\left(\left\langle y_l\left(t, s, \omega, 
\theta\right), \sqrt{2\lambda}\cdot\right\rangle\right)(\theta)\\
 & =d_{\theta}\left\langle y_l(t, s, \theta, \omega),
\sqrt{2\lambda}\cdot\right\rangle(\theta).
\end{aligned}
\end{equation}
From (\ref{somega}) and (\ref{stheta}) we therefore have
that $S_l$ is
a non-degenerate phase function such that $SR_l(\lambda)=\Lambda_{S_{l}}.$ $\hfill\Box$

\subsection{Resolvent Relation}\label{sresolrel}
We now define the Lagrangian submanifold which we will prove is quantized 
by the cut-off
resolvent in the sense of semi-classical Fourier integral operator.

We set
\begin{equation*}
\tilde{\Lambda}_{R}\left(\lambda, J\right)=\cup_{t\in
J}\graph\exp\left(tH_p\right)|_{\Sigma_{\lambda}},
\end{equation*}
where $J\subset\mathbb{R}$ is an open interval.
We assume that for every $\left(x, \xi\right)\in\Sigma_{\lambda}$ and
every $t\in~J,$
$\exp \left(tH_{p}\right)\left(\left(x, \xi\right)\right)\ne \left(x,
\xi\right).$
Then $\tilde{\Lambda}_{R}\left(\lambda, J\right)$ is a
Lagrangian
submanifold
of
$\left(T^{*}\mathbb{R}^{n}\times T^{*}\mathbb{R}^{n},
\tilde{\sigma}\right),$
where  $\tilde{\sigma}=\pi_{1}^{*}\sigma-\pi_{2}^{*}\sigma=d\xi\wedge dx -
d\eta\wedge
dy,$ where $\pi_{j}: T^{*}\mathbb{R}^{n}\times T^{*}\mathbb{R}^{n}\to
T^{*}\mathbb{R}^{n},$ $j=1, 2$ is the projection onto the $j-$th factor.
We define the resolvent relation as
$\Lambda_{R}(\lambda)=\tilde{\Lambda}_{R}'\left(\lambda,
\mathbb{R}^{+}\right)\cap T^{*}\left(\supp \tilde{\chi}_{2}\times
\supp
\tilde{\chi}_{1}\right).$

\section{Proof of Main Theorem}\label{cprm}

We now turn to the proof of our Main Theorem.

\medskip
\noindent
{\it Proof:}\;
We first prove that the scattering amplitude belongs to 
$\mathcal{D}_{h}'(\mathbb{S}^{n-1}\times\mathbb{S}^{n-1}).$
For that, let $\psi\in C^{\infty}(\mathbb{S}^{n-1}\times\mathbb{S}^{n-1})$ 
have support in a coordinate chart on 
$\mathbb{S}^{n-1}\times\mathbb{S}^{n-1}$ with local coordinates 
$\left(\tilde{\omega}, \tilde{\theta}\right)$ and let $\psi_{1},$ 
$\psi_{2}\in C^{\infty}(\mathbb{S}^{n-1})$ be such that $\psi_{1}\times 
\psi_{2}=1$ on 
$\supp\psi.$
Then, since $K_{A(\lambda, h)}\in 
C^{\infty}(\mathbb{S}^{n-1}\times\mathbb{S}^{n-1}),$ we have
\begin{equation*}
\begin{aligned}
\left|\mathcal{F}_{h}\left(K_{A(\lambda, h)}\right)\left(\tilde{\xi}, 
\tilde{\eta}\right)\right|& = \left|\int\int K_{A(\lambda, 
h)}\left(\tilde{\omega}, \tilde{\theta}\right)\psi\left(\tilde{\omega}, \tilde{\theta}\right)e^{-\frac{i}{h}
\left(\langle \tilde{\omega}, \tilde{\xi}\rangle+\langle \tilde{\theta}, 
\tilde{\eta}\rangle\right)} d\tilde{\theta} d\tilde{\omega}\right|\\
& \leq\int\int\left| K_{A(\lambda, h)}\left(\tilde{\omega},
\tilde{\theta}\right)\psi\left(\tilde{\omega}, \tilde{\theta}\right)\right|d\tilde{\theta} d\tilde{\omega}\\
& \leq C \int\int\left| K_{A(\lambda, h)}\left(\tilde{\omega},
\tilde{\theta}\right)\psi_{1}\left(\tilde{\omega}\right)\psi_{2}\left(\tilde{\theta}\right)\right|d\tilde{\theta} d\tilde{\omega}\\
& \leq C \left\|A(\lambda, 
h)\right\|_{\mathcal{B}(L^{2}(\mathbb{S}^{n-1}))}\left\|\psi_{1}\right\|_{L^{2}(\mathbb{S}^{n-1})}
\left\|\psi_{2}\right\|_{L^{2}(\mathbb{S}^{n-1})}\\
& =\mathcal{O}\left(h^{\frac{1-n}{2}}\right),
\end{aligned}
\end{equation*} 
where the last equality follows from 
\[A(\lambda, h)=h^{\frac{n-1}{2}}\lambda^{\frac{1-n}{4}}2^{\frac{3n-13}{4}}
\pi^{\frac{n-3}{2}}e^{i\pi\frac{5-n}{4}}(S(\lambda, h)-I)\]
 and the fact that 
$S(\lambda, h)$ is a unitary operator on $L^{2}(\mathbb{S}^{n-1}).$ 

We also observe that a direct calculation shows that 
$\pi_{1}\circ\left(\pi_{2}|_{\Lambda(\lambda)}\right)^{-1}
\left(\Lambda_{R}\left(\lambda\right)\right)$ is a Lagrangian submanifold of 
$T^{*}\mathbb{S}^{n-1}\times T^{*}\mathbb{S}^{n-1}.$
To prove the theorem, now, we first note that 
$$\mathbb{E}_{-}(\lambda, 
h)\otimes\mathbb{E}_{+}(\lambda, h)\in 
\mathcal{I}_{h}^{-n+\frac{1}{2}}(\mathbb{S}^{n-1}\times 
\mathbb{S}^{n-1}\times\mathbb{R}^{n}\times \mathbb{R}^{n}, 
\Lambda(\lambda))$$ and we easily see that 
\begin{equation}\label{immerse}
\pi_{2}|_{\Lambda(\lambda)}\text{ is an immersion.}
\end{equation}
Let, now,   
\[A_{j}\in\Psi_{h}^{0}\left(1, 
\mathbb{S}^{n-1}\times\mathbb{S}^{n-1}\right),\;
j=1, \dots, N\] have symbols supported in a neighborhood of $\bar{p}$ and principal symbols vanishing on 
$\pi_{1}\circ\left(\pi_{2}|_{\Lambda(\lambda)}\right)^{-1}
\left(\Lambda_{R}\left(\lambda\right)\right).$
From \eqref{immerse} and \cite[Lemma 7]{AI} we deduce that there exist 
$D, B_{j}\in\Psi_{h}^{0}\left(1, \mathbb{R}^{2n}\right),$ $j=1, \dots, N,$ with 
symbols supported near the point 
\[\bar{q}=\left(x_0+t\omega_0, y_0+s\theta_0, \sqrt{2\lambda}\omega_0, -\sqrt{2\lambda}\theta_0 \right)\text{ for some } s\in[s_1, s_2],   t\in[t_1, t_2],\]  such that 
\begin{equation}\label{symbB}
\left(\pi_{2}|_{\Lambda(\lambda)}\right)^{*}\sigma_{0}\left(B_{j}\right)
=-\left(\pi_{1}|_{\Lambda(\lambda)}\right)^{*}\sigma_{0}\left(A_{j}\right), \; j=1, \dots, N,
\end{equation}
\begin{equation*}
\left(\pi_{2}|_{\Lambda(\lambda)}\right)^{*}\sigma_{0}\left(D\right)
=-\left(\pi_{1}|_{\Lambda(\lambda)}\right)^{*}\sigma_{0}\left(C\right),
\end{equation*}
and
\begin{equation}\label{aeeb}
\left(\prod_{j=1}^{N} A_{j}\right)C
\left(\mathbb{E}_{-}\left(\lambda, h\right)\otimes
\mathbb{E}_{+}\left(\lambda, h\right)\right)\chi
\equiv\left(\mathbb{E}_{-}\left(\lambda, h\right)\otimes 
\mathbb{E}_{+}\left(\lambda, h\right)\right)\chi\cdot
\left(\prod_{j=1}^{N} B_{j}\right)D
\end{equation}
near $\left(\bar{p}, \bar{q}\right),$
where $\chi\in C_{c}^{\infty}(\mathbb{R}^{n}\times\mathbb{R}^{n})$ is equal to 1 on 
$\supp\tilde{\chi}_{2}\times\supp\tilde{\chi}_{1}.$

Let, now, $Y=\left[h^{2}\Delta, \chi_{2}\right]\otimes \left[h^{2}\Delta,
\chi_{1}\right]^{t}.$
Since $Y\in\Psi_{h}^{-1, 1}(\mathbb{R}^{2n})$ we obtain from \cite[Lemma 5]{AI} that
\begin{equation}\label{yr}
YK_{\tilde{\chi}_2 R\left(\lambda, h\right)\tilde{\chi}_1}\in 
I_{h}^{r-1}(\mathbb{R}^{n}\times\mathbb{R}^{n}, \Lambda_{R}(\lambda)).
\end{equation}
Therefore, from \eqref{aeeb}, we have
\begin{equation}\label{maininterwine}
\begin{aligned}
& \left(\prod_{j=1}^{N} A_{j}\right)C
K_{A\left(\lambda, h\right)}\\
& \quad=c\left(n, \lambda, h\right)
\left(\mathbb{E}_{-}\left(\lambda, h\right)\otimes 
\mathbb{E}_{+}\left(\lambda, h\right)\right)\chi\cdot
\left(\prod_{j=1}^{N} B_{j}\right) D
YK_{\tilde{\chi}_2 R\left(\lambda, h\right)\tilde{\chi}_1}+\mathcal{O}_{L^{2}(\mathbb{S}^{n-1}\times\mathbb{S}^{n-1})}(h^{\infty}).
\end{aligned}
\end{equation}
The choice of the operators $A_{j}$ and  \eqref{symbB} now imply 
\begin{equation*}
\sigma_{0}\left(B_{j}\right)|_{\Lambda_{R}\left(\lambda\right)}=0, j=1, 
\dots, N.
\end{equation*}
From \cite[Lemma 5]{AI} we have again 
\[DYK_{\tilde{\chi}_2 R\left(\lambda, 
h\right)\tilde{\chi}_1}\in I_{h}^{r-1}(\mathbb{R}^{n}\times\mathbb{R}^{n}, \Lambda_{R}(\lambda)).\]
Therefore,
\begin{equation}\label{presub}
\left(\prod_{j=1}^{N} B_{j}\right)DYK_{\tilde{\chi}_2 R\left(\lambda,
h\right)\tilde{\chi}_1}=\mathcal{O}_{L^{2}(\mathbb{R}^{2n})}\left(h^{N-r+1+\frac{n}{4}}\right).
\end{equation}

Lastly, as in the proof of \cite[Proposition 2.1]{Burq}, we have that 
\begin{equation}\label{normE}
\left\|\mathbb{E}^{\phi}_{\pm}\right\|_{\mathcal{B}(L^{2}(\mathbb{R}^{n}), 
L^{2}(\mathbb{S}^{n-1}))}=\mathcal{O}\left(h^{\frac{n-1}{2}}\right),
\end{equation}
where $\mathbb{E}_{\pm}^{\phi}\left(\lambda, h\right)$ are the operators 
with Schwartz kernels \[K_{\mathbb{E}_{\pm}^{\phi}\left(\lambda, 
h\right)}(\omega, x)=\phi\left(x\right)\exp
\left(\pm i\sqrt{2\lambda}\left\langle \omega, x \right\rangle/h\right)\] 
for $\phi\in 
C_c^{\infty}\left(\mathbb{R}^{n}\right).$

We now substitute (\ref{presub}) into (\ref{maininterwine}) and use 
(\ref{normE}) to obtain 
\begin{equation*}
\left(\prod_{j=1}^{N} A_{j}\right)CK_{A\left(\lambda, h\right)}= 
\mathcal{O}_{L^{2}\left(\mathbb{R}^{2\left(n-1\right)}\right)}\left(h^{N-r-\frac{n}{2}}\right).
\end{equation*}
Therefore $CK_{A(\lambda, h)}\in 
I_{h}^{r+\frac{1}{2}}\left(\mathbb{S}^{n-1}\times\mathbb{S}^{n-1}, 
\pi_{1}\left(\pi_{2}|_{\Lambda(\lambda)}\right)^{-1}(\Lambda_{R}(\lambda))\right).$ 

\section{Applications}\label{cappl}
Here we discuss two applications of the main theorem to compactly
supported potential and metric perturbations of the Euclidean Laplacian.
The setting is as follows.

Let the setting be as in Section \ref{scr}.
The operators 
$P\left(h\right)=\frac{1}{2}h^{2}\Delta_{g}+V,$ $0<h\leq 1,$ acting on 
$\mathcal{H}=L^{2}\left(\mathbb{R}^{n}, d\text{vol}_{g}\right)$ 
and equipped with the 
common domain
$\mathcal{D}=H^{2}\left(\mathbb{R}^{n}, d \text{vol}_{g}\right)$ 
admit unique self-adjoint extensions, which we 
denote by the
same notation.
As before, we denote their resolvents by $R\left(z, 
h\right)=\left(P(h)-z\right)^{-1}$ for
$z\in\mathbb{C}_+$ and use the same notation for their meromorphic
continuations.
Proposition 2.3 in \cite{Mel} states that there are no resonances in 
$\mathbb{R}\backslash\{0\}$ in the case of a smooth compactly supported potential 
perturbation.

We define a non-trapping energy level as follows:
\begin{Def}
Let $\left\{\left(x\left(\cdot; x_{0}, \xi_{0}\right), \xi\left(\cdot; 
x_{0}, 
\xi_{0}\right)\right)\right\}$ be the
integral curve of
$H_{p}$ with initial conditions $\left(x_{0}, \xi_{0}\right).$
The energy $\lambda>0$ is non-trapping if for every $r>0$, there exists
$s>0$ such that $\left(x_{0}, 
\xi_{0}\right)\in\Sigma_{\lambda}$ with
$\|x_{0}\|<r$ implies that
$\|x\left(s; x_{0}, \xi_{0}\right)\|>r$ for every $|s|>t.$
We also introduce the notation $T(r)$ for the infimum over $s$ 
with this property.
\end{Def}

\subsection{The Cut-off Resolvent as a Semi-Classical Fourier Integral 
Operator}\label{sresolv}

We now prove that the second assumption of our main theorem is satisfied 
in the setting we have just described.

\begin{Th}\label{resolventfio}
Let $\|\chi R\left(\lambda, h\right)\chi\|_{\mathcal{B}\left(L^{2}\left(\mathbb{R}^{n}\right)\right)}\leq 
C_{\chi} h^{s}$ for every $\chi\in 
C^{\infty}_{c}\left(\mathbb{R}^{n}\right)$ and 
some 
$s\in\mathbb{R}.$ 
Let $\rho_{0}\in\Lambda_{R}(\lambda)$ be such that $\gamma(\cdot; 
\pi_{1}(\rho_{0}))$ 
is a non-trapped trajectory.

Then there exists an open set 
$V\subset\Lambda_{R}(\lambda),$ $\rho_{0}\in V,$ such that 
$$\tilde{\chi}_{2} 
R\left(\lambda, h\right)\tilde{\chi}_{1}\in 
\mathcal{I}_{h}^{1}\left(\mathbb{R}^{2n}, \overline{V}\cap~\Lambda_{R}(\lambda)\right).$$ 
\end{Th}

\medskip
\noindent
{\it Proof:}\; 
First, we prove that $K_{\tilde{\chi}_{2} R\left(\lambda, h\right)\tilde{\chi}_{1}}\in\mathcal{D}_{h}'(\mathbb{R}^{2n}).$
For this, let $\chi\in C_{c}^{\infty}\left(\mathbb{R}^{2n}\right)$ and let $\rho_{j}\in 
C_{c}^{\infty}(\mathbb{R}^{n}\backslash B(0, R_{0})),$ $j=1, 2,$ be such that 
$\rho_2\times\rho_1=1$ on 
$\supp\tilde{\chi}_{2}\times\supp\tilde{\chi}_{1}\cap\supp\chi.$
Then, since $K_{\tilde{\chi}_{2}R\left(\lambda, h\right)\tilde{\chi}_{1}}\in 
C^{\infty}(\mathbb{R}^{2n}),$ we have
\begin{equation*} 
\begin{aligned}
\left|\mathcal{F}_{h}\left(K_{\tilde{\chi}_{2}R\left(\lambda, 
h\right)\tilde{\chi}_{1}}\right)(\xi, \eta) \right|&= 
\left|\int\int K_{\tilde{\chi}_{2}R\left(\lambda, h\right)\tilde{\chi}_{1}}(x, y)\chi(x, 
y)e^{-\frac{i}{h}\left(\langle x, \xi\rangle+\langle y, \eta \rangle\right)}dxdy \right|\\
& \leq \int\int\left|K_{\tilde{\chi}_{2}R\left(\lambda, h\right)\tilde{\chi}_{1}}(x, 
y)\chi(x, y)\right|dxdy\\
& \leq \int\int\left|K_{\rho_{2}R\left(\lambda, 
h\right)\rho_{1}}(x, y)\left(\rho_{2}\otimes\rho_{1}\right)(x, y)\right| d 
x d y\\
& \leq 
\left\|\rho_{2}\right\|_{L^{2}(\mathbb{R}^{n})}\left\|\rho_{1}\right\|_{L^{2}(\mathbb{R}^{n})}
\left\|\rho_{2}R\left(\lambda, 
h\right)\rho_{1}\right\|_{\mathcal{B}(L^{2}(\mathbb{R}^{n}))}\\
&=\mathcal{O}(h^{s}),
\end{aligned}
\end{equation*}
which verifies the assertion.

To recall the representation of the resolvent, which we shall use to prove the lemma, we 
recall that for $f\in 
C_{c}^{\infty}(\mathbb{R}^{n})$
\begin{equation*}
-\frac{i}{h}\left(P\left(h\right)-\lambda\right)\int_{0}^{T}e^{\frac{i}{h}t\lambda}U(t)f
dt
=e^{\frac{i}{h}T\lambda} U(T)f-f,
\end{equation*}
where $U(t)=e^{-\frac{i}{h}tP(h)},$ $t\in\mathbb{R}$ is the unitary group of $P(h).$
The same proof as in Lemma B.1, \cite{Michel2}, shows that 
$\left(1-\chi_{0}\right)U(T)f\in\mathcal{S}(\mathbb{R}^{n}),$ $\chi_0\in 
C_{c}^{\infty}(\mathbb{R}^{n}),$ $\chi_{0}=1$ on $B(0, R_{0}).$
Since we can also think of $R(\lambda, h)$ as the limit $\lim_{\epsilon\to 
0, \epsilon>0}R(\lambda\pm i\epsilon, h)$ in the spaces of bounded 
operators $\mathcal{B}\left(L^{2}_{\alpha}, L^{2}_{-\alpha}\right),$ 
$\alpha>\frac{1}{2},$ where 
\[L^{2}_{\alpha}=\left\{f: \left(\chi_{0}+(1-\chi_{0})\langle 
x\rangle^{\alpha}\right)f\in L^{2}(X, 
d\text{vol}_{g}\mathbb{R}^{n})\right\},\] we obtain 
\begin{equation}\label{resrepr} 
\tilde{\chi}_{2}R\left(\lambda, 
h\right)\tilde{\chi}_{1}=\frac{i}{h}\int_{0}^{T}e^{\frac{i}{h}t\lambda}\tilde{\chi}_{2}
U(t)\tilde{\chi}_{1}dt+e^{\frac{i}{h}T\lambda}\tilde{\chi}_{2}R\left(\lambda,
h\right)U(T)\tilde{\chi}_{1}
\end{equation}
and this is the representation of the resolvent, which we shall use in 
this proof.

We further recall the well-known fact that $U(t)\in \mathcal{I}_{h}^{0}(\mathbb{R}^{2n}, 
\Lambda_{t}),$ $t\in\mathbb{R},$ where $\Lambda_t=\graph\exp(t H_p).$
Since $\rho_{0}\in\Lambda_{R}\left(\lambda\right)$ is such that $\gamma\left(\cdot; 
\pi_{1}(\rho_0)\right)$ is a non-trapped trajectory, there exists an open 
set 
$V\subset\Lambda_{R}\left(\lambda\right),$ $\rho_0\in V,$ such that for 
every $\rho\in V,$ 
$\gamma\left(\cdot; \pi_{1}(\rho)\right)$ is a non-trapped trajectory. 
By adjusting $V,$ if necessary, we can assume that the same holds for 
all points in $\bar{V}.$
Let $Q\in\Psi_{h}^{0}\left(1, 
\mathbb{R}^{2n}\right)$ be a microlocal cut-off to the
neighborhood $V$ as in (\ref{P}) with a compactly 
supported symbol.
First, we shall prove that $QK_{\tilde{\chi}_{2}R\left(\lambda, 
h\right)U(T)\tilde{\chi}_{1}}=\mathcal{O}_{L^{2}(\mathbb{R}^{2n})}(h^{\infty})$ 
for $T>0$ sufficiently large.
By \cite[Lemma 3, (a)]{AI} and the choice of $Q,$ we have that 
\begin{equation}\label{remwfi}
WF_{h}^{i}\left(QK_{\tilde{\chi}_{2}R\left(\lambda,
h\right)U(T)\tilde{\chi}_{1}}\right)=\emptyset,
\end{equation}
and therefore by Proposition 
7.1, (i), \cite{Michel}, it follows that it is sufficient to prove that there 
exists $\bar{T}>0$ such that for every
$T>\bar{T},$ $WF_{h}^{f}\left(QK_{\tilde{\chi}_{2}R\left(\lambda, h\right)
e^{-\frac{i}{h}TP\left(h\right)}\tilde{\chi}_{1}}\right)=\emptyset.$
This will follow, if we prove that 
$WF_{h}^{f}\left(K_{\tilde{\chi}_{2}R\left(\lambda, h\right)
e^{-\frac{i}{h}TP\left(h\right)}\tilde{\chi}_{1}}\right)\cap 
WF_{h}^{f}(Q)=\emptyset.$
To prove the latter, consider 
\begin{equation*}
\left\langle K_{\tilde{\chi}_{2}R\left(\lambda, h\right)
e^{-\frac{i}{h}TP\left(h\right)}\tilde{\chi}_{1}}, 
\left(\psi_{2}\otimes\psi_{1}\right)e^{-\frac{i}{h}\left(\langle \cdot, 
\xi\rangle+\langle \cdot\cdot, \eta\rangle\right)}\right\rangle, 
\end{equation*} 
where $\supp\psi_{2}\times \{\xi\}\times\supp\psi_{1}\times\{\eta\}\subset 
V.$
Let $V_{1}\subset\mathbb{R}^{n}$ be a bounded 
open set such that $\supp\psi_{2}\times \{\xi\}\times\supp\psi_{1}\times 
V_{1}\subset V,$ $\eta\in V_1.$
Now for every $\eta$ we have that 
$WF_{h}\left(\psi_{1}e^{-\frac{i}{h}\langle\cdot, 
\eta\rangle}\right)=\supp\psi_{1}\times\{\eta\}$ is 
compact.
This, together with \cite[Theorem 1]{AI}, \cite[Lemma 6]{AI}, and \cite[Lemma 3, (c)]{AI}, and the fact that $U(t)\in \mathcal{I}_{h}^{0}(\mathbb{R}^{2n}, \Lambda_t),$ $t\in\mathbb{R},$ allows us to conclude that
\begin{equation}\label{propu(t)}
WF_{h}^{f}\left(U(T)\tilde{\chi}_{1}\psi_{1}e^{-\frac{i}{h}\langle\cdot, 
\eta\rangle}\right)\subset \exp\left(T 
H_{p}\right)\left(WF_{h}^{f}\left(\psi_{1}e^{-\frac{i}{h}\langle\cdot, 
\eta\rangle}\right)\right), \eta\in V_1. 
\end{equation} 
After decreasing $V_{1},$ if necessary, we have, by the proof of\cite[Lemma 4]{AI}, 
that the estimates in (\ref{propu(t)}) can be made uniform in $\eta\in 
V_{1}.$
Since $\exp\left(T H_{p}\right)$ is a diffeomorphism, it follows
that 
$\cup_{\eta\in\bar{V}_{1}} 
\exp\left(T 
H_{p}\right)\left(WF_{h}^{f}\left(\psi_{1}e^{-\frac{i}{h}\langle\cdot,
\eta\rangle}\right)\right)$ is compact.
Further, as we are working with the outgoing resolvent, we have that
\begin{equation*}
\begin{aligned}
WF_{h}^{f}\left(R(\lambda, h)U(T)\tilde{\chi}_{1}\psi_{1}e^{-\frac{i}{h}\langle\cdot,
\eta\rangle}\right) & 
\subset\cup_{t>0}\exp\left(tH_{p}\right)\left(WF_{h}^{f}\left(U(T) \tilde{\chi}_{1}\psi_{1}e^{-\frac{i}{h}\langle\cdot, 
\eta\rangle}\right)\right)\\
 & \subset\cup_{t>0}\exp\left((t+T)H_{p}\right)\left(WF_{h}^{f}\left(\psi_{1}
e^{-\frac{i}{h}\langle\cdot, \eta\rangle}\right)\right). 
\end{aligned}
\end{equation*}
By the non-trapping assumption, there exists $\bar{T}>0$ such that for 
every $T>\bar{T}$ and every $(y, \eta)\in \supp\psi_{1}\times V_{1}$ we have $x(T; y, 
\eta)\in \left(\supp\tilde{\chi}_{2}\right)^{c}.$ 
We now let $T>\bar{T}$ be fixed and we have
\begin{equation*}
\left\langle K_{\tilde{\chi}_{2}R\left(\lambda, h\right)
e^{-\frac{i}{h}TP\left(h\right)}\tilde{\chi}_{1}},
\left(\psi_{2}\otimes\psi_{1}\right)e^{-\frac{i}{h}\left(\langle \cdot,
\xi\rangle+\langle \cdot\cdot, \eta\rangle\right)}\right\rangle=\mathcal{O}(h^{\infty})
\end{equation*}
for every $\eta\in V_1$ and uniformly in $\xi\in U,$ where
$U\subset\mathbb{R}^{n}$ a bounded open set such that $\supp\psi_{2}\times 
U\times\supp\psi_{1}\times V_{1}\subset V.$
The proof of \cite[Lemma 4]{AI}, now shows again that the estimate here 
can be made uniform 
in $\eta\in V_{1}.$
We thus have that for every $T>\bar{T}$ 
\begin{equation}\label{nowfh}
WF_{h}^{f}\left(QK_{\tilde{\chi}_{2}R(\lambda, 
h)U(T)\tilde{\chi}_{1}}\right)=\emptyset,
\end{equation}
which, together with (\ref{remwfi}), gives
\begin{equation}\label{remainder}
QK_{\tilde{\chi}_{2}R(\lambda, h)U(T)\tilde{\chi}_{1}}=\mathcal{O}_{L^{2}\left(\mathbb{R}^{2n}\right)}(h^{\infty}).
\end{equation}

Let, now, $\tilde{Q}\in\Psi_{h}^{0}(1, \mathbb{R}^{2n})$ have a compactly 
supported symbol.
Then, since $U(t)\in\mathcal{I}_{h}^{0}(\mathbb{R}^{2n}, \Lambda_{t}),$ $t\in\mathbb{R},$ we have 
\begin{equation*}
\tilde{Q}QK_{\tilde{\chi}_{2}R(\lambda, h)\tilde{\chi}_{1}}=\frac{i}{h}\int_{0}^{T}e^{\frac{i}{h}t\lambda}\left(\tilde{\chi}_{2}\otimes
\tilde{\chi}_{1}\right)K_{U(t)}dt+e^{\frac{i}{h}T\lambda}
\tilde{Q}QK_{\tilde{\chi}_{2}R(\lambda, h)U(T)\tilde{\chi}_{1}}= \mathcal{O}_{L^{2}(\mathbb{R}^{2n})}\left(h^{-1-\frac{n}{2}}\right).
\end{equation*}

Let, now, $B_{j}\in\Psi_{h}^{0}\left(1, \mathbb{R}^{2n}\right),$ $j=1, 
\dots, 
k$ have compactly supported symbols and principal symbols 
$b_{j}$ 
vanishing on $\Lambda_{R}\left(\lambda\right)$ and consider
\begin{equation}\label{bs}
\left(\prod_{j=1}^{k} B_{j}\right)QK_{\tilde{\chi}_{2}R(\lambda, h)\tilde{\chi}_{1}}.
\end{equation}
We have that 
$\Lambda_{R}\left(\lambda\right)=\left(\cup_{s\in\mathbb{R}}\Lambda_{s}\right)\cap 
\left(\Sigma_{\lambda}\times \Sigma_{\lambda}\right),$ and the intersection at every 
point is clean.
Therefore, by Proposition C.3.1 \cite{H}, vol.~3, we can choose local 
coordinates around 
$\rho$ such that $\cup_{s\in\mathbb{R}}\Lambda_{s}$ and 
$\Sigma_{\lambda}\times 
\Sigma_{\lambda}$ are given there by linear equations in the local 
coordinates.
This implies that for every $j=1, \dots, k$ we can find functions $c_j, 
g_j, h_j\in 
C^{\infty}_{c}\left(T^{*}\mathbb{R}^{n}\right)$ with $c_j$ vanishing on 
$\cup_{s\in\mathbb{R}}\Lambda_{s},$ 
$j=1, \dots, k,$ such that $b_j=c_j+\left(\left(p-\lambda\right)\otimes 
1\right)g_j+\left(1\otimes\left(p-\lambda\right)\right)h_j.$
Now, for $a, b\in S^{0}_{2n}(1)$ we have that 
$Op_{h}(a)Op_{h}(b)=Op_{h}(ab)+\mathcal{O}_{\mathcal{B}(L^{2}(\mathbb{R}^{n}))}(h)$ 
and we can 
therefore rewrite (\ref{bs}) as follows
\begin{equation*}
\begin{aligned}
\left(\prod_{j=1}^{k} B_{j}\right)QK_{\tilde{\chi}_{2}R(\lambda, 
h)\tilde{\chi}_{1}} & =\Bigg(\prod_{j=1}^{k} (Op_{h}(c_j)+Op_{h}(g_j)((P(h)-\lambda)\otimes I) \\
& \quad + Op_{h}(h_j)(I\otimes(P(h)-\lambda))+hS_j)\Bigg)QK_{\tilde{\chi}_{2}R(\lambda, 
h)\tilde{\chi}_{1}},
\end{aligned} 
\end{equation*}
where $S_{j}\in\Psi_{h}^{0}(1, \mathbb{R}^{2n})$ and $\sigma(S_j)\in 
C_{c}^{\infty}(\mathbb{R}^{4n}),$ $j=1, \dots, k.$
This we further rewrite as
\begin{equation}\label{prods}
\left(\prod_{j=1}^{k} B_{j}\right)QK_{\tilde{\chi}_{2}R(\lambda,
h)\tilde{\chi}_{1}}=\sum_{l=1}^{4^{k}} \left(\prod_{j=1}^{k} 
T^{l}_{j}\right)QK_{\tilde{\chi}_{2}R(\lambda,
h)\tilde{\chi}_{1}},
\end{equation}
where 
\[T^{l}_{j}\in\left\{Op_{h}(c_j), Op_{h}(g_j)((P(h)-\lambda)\otimes
I), Op_{h}(h_j)(I\otimes(P(h)-\lambda)), h S_j\right\},\] 
$l=1, \dots, 
4^k,$ 
$j=1, \dots, k.$
We now turn to analyzing the individual summands.
As the superscript will not be 
important, we will omit it from the notation.

We consider the case $k>1.$
The case $k=1$ will be implicit in the discussion below.
Let $m\in\{2, \dots, k\}$ be the 
largest index such that 
$$T_{m}\in\{Op_{h}(c_{m}), h S_{m}\}$$ and 
$$T_{m-1}\in\{Op_{h}(g_{m-1})((P(h)-\lambda)\otimes
I), Op_{h}(h_{m-1})(I\otimes(P(h)-\lambda))\}.$$ 
We first assume that $T_{m}=Op_{h}(c_{m}).$
Since $$\sigma_{0}([T_m, T_{m-1}])=\frac{h}{i}\{\sigma_{0}(T_m), 
\sigma_0(T_{m-1})\}$$ and $\sigma_{0}(T_m)$ and $\sigma_0(T_{m-1})$ 
vanish on the Lagrangian submanifold $\Lambda_{R}(\lambda)$ near $\rho,$ it follows that $\sigma_{0}([T_m, T_{m-1}])$ also vanishes 
on 
$\Lambda_{R}(\lambda)$ near $\rho.$
Therefore we have, as before, that $$\sigma_{0}([T_m, T_{m-1}])
=\tilde{c}+\tilde{g}((p-\lambda)\otimes 1)+\tilde{h}(1\otimes 
(p-\lambda))$$ for some $\tilde{c}, \tilde{g}, \tilde{h}\in 
C_{c}^{\infty}(\mathbb{R}^{4n})$ with supports in a sufficiently small 
neighborhood of $\rho$ and 
$\tilde{c}|_{\cup_{s\in\mathbb{R}}\Lambda_{s}}=0.$ 
Thus 
\begin{equation}\label{subprod1}
\begin{aligned}
& T_{m-1}T_m\\
& \quad =Op_{h}(\tilde{c})+Op_{h}\left(\tilde{g}\right)((P(h)-\lambda)\otimes 
I)+Op_{h}\left(\tilde{h}\right)(I\otimes (P(h)-\lambda))\\
& \quad\quad +T_m T_{m-1}+h\tilde{T}
\end{aligned}
\end{equation}
where $\tilde{T}\in \Psi_{h}^{0}(1, \mathbb{R}^{2n})$ and
$\sigma\left(\tilde{T}\right)\in C_{c}^{\infty}(\mathbb{R}^{4n})$ and we
use the latter expression in (\ref{subprod1}) to replace $T_{m-1}T_m$ in
the product above.

If, now, $T_{m}=h S_{m}$ we rewrite $T_{m-1}T_m=T_{m}T_{m-1}+[T_{m-1}, 
T_m]$ in the product
above and observe that $[T_{m-1}, T_m]\in\Psi_{h}^{-2}(1, 
\mathbb{R}^{4n}).$

We iterate this process until each product which appears in (\ref{prods}),
where we may now have more than $4^{k}$ products, is of the form
\begin{equation*}
h^{k-k_{1}}\prod_{j=1}^{k_{1}}Q_j,\;  k_{1}\in\{1, \dots, k\},
\end{equation*}
where for some $j_{0}\in\{0, \dots, k_{1}\}$ we have that
\begin{equation*}
\begin{aligned}
& \text{for } j_0<j\leq k_{1},\, Q_{j} \in 
\left\{Op_{h}\left(g^{1}_{j}\right)
((P(h)-\lambda)\otimes I), Op_{h}\left(h^{1}_{j}\right)(I\otimes 
(P(h)-\lambda))\right\},\\
& \quad \text{ for some } g^{1}_{j}, h^{1}_{j}\in 
C^{\infty}_{c}(\mathbb{R}^{4n}),\\
& \text{for } 1\leq j\leq j_0, \, Q_{j} \in \left\{Op_{h}(c^{1}_j), h 
S_{j}^{1}\right\}, \\
& \quad \text{for some } c^{1}_{j}\in C_{c}^{\infty}(T^{*}\mathbb{R}^{n}),
c^{1}_{j}|_{\cup_{t\in\mathbb{R}}\Lambda_{t}}=0, 
S_{j}^{1}\in\Psi_{h}^{0}(1, \mathbb{R}^{2n}),
\sigma( S_{j}^{1})\in C_{c}^{\infty}(\mathbb{R}^{4n}),
\end{aligned}
\end{equation*}
with all symbols supported in a sufficiently small neighborhood of $\rho.$

Next, we let $\tilde{m}\in\{2, \dots, k_{1}\}$ denote the largest
index for which we have $T_{\tilde{m}-1}=Op_{h}(c^{1}_{\tilde{m}-1})$ and
$T_{\tilde{m}}=h S_{\tilde{m}}^{1}.$
We then replace $T_{\tilde{m}-1}T_{\tilde{m}}$ by 
$T_{\tilde{m}}T_{\tilde{m}-1}
+[T_{\tilde{m}-1}, T_{\tilde{m}}]$ in the product above and
observe that $[T_{\tilde{m}-1}, T_{\tilde{m}}]\in \Psi_{h}^{-1}(1,
\mathbb{R}^{2n})$ and
$\sigma\left( [T_{\tilde{m}-1}, T_{\tilde{m}}]\right)\in
C_{c}^{\infty}\left(\mathbb{R}^{4n} \right).$

We repeat this procedure until every product which appears in 
(\ref{prods}) is of the form
\begin{equation*}
h^{k-k_{2}} \prod_{j=1}^{k_{2}}V_{j}, \; k_{2}\in\{1, \dots, k\},
\end{equation*}
where for some $j_{1}, j_{2}\in\{0, \dots, k_{2}\},$ $j_{1}\leq
j_{2},$ we have that
\begin{equation*}
\begin{aligned}
& \text{for } j_{2}<j\leq k_{2}, V_{j} \in 
\left\{Op_{h}\left(g^{2}_{j}\right)
((P(h)-\lambda)\otimes I),
Op_{h}\left(h^{2}_{j}\right)(I\otimes (P(h)-\lambda))\right\},\\
& \quad \text{for some } g^{2}_{j}, h^{2}_{j}\in 
C^{\infty}_{c}(\mathbb{R}^{4n}),\\
& \text{for } j_{1}<j\leq j_{2}, V_{j} = Op_{h}(c^{2}_j),
\text{ for some } c^{2}_{j}\in C_{c}^{\infty}(T^{*}\mathbb{R}^{n}),
c^{2}_{j}|_{\cup_{t\in\mathbb{R}}\Lambda_{t}}=0,\\
& \text{for } 1\leq j\leq j_{1}, V_{j} = h S_{j}^{2},
\text{ for some } S_{j}^{2}\in\Psi_{h}^{0}(1, \mathbb{R}^{2n}),
\sigma(S_{j}^{2})\in C_{c}^{\infty}(\mathbb{R}^{4n}),
\end{aligned}
\end{equation*}
where again all symbols are supported in a sufficiently small neighborhood 
of $\rho.$

We shall again omit the superscripts from the notation below.
We also observe that the symmetry of $K_{\tilde{\chi}_{2}R\left(\lambda,
h\right)\tilde{\chi}_{1}}$ allows us to assume that
\begin{equation*}
V_{j}=Op_{h}\left(g_{j}\right)((P(h)-\lambda)\otimes I), \; g_{j}\in
C^{\infty}_{c}(\mathbb{R}^{4n}), \; j_{2}<j\leq k_{2}.
\end{equation*}

We now analyze
\begin{equation}\label{1pl}
\begin{aligned}
(\left(P(h)-\lambda\right) & \otimes
I)QK_{\tilde{\chi}_{2}R\left(\lambda,
h\right)\tilde{\chi}_{1}}\\
& =\left[\left(P(h)-\lambda\right)\otimes I,
Q\right]K_{\tilde{\chi}_{2}R\left(\lambda,
h\right)\tilde{\chi}_{1}}+Q\left(\left(-h^{2}\Delta-\lambda\right)\otimes
I\right)K_{\tilde{\chi}_{2}R\left(\lambda,
h\right)\tilde{\chi}_{1}}
\end{aligned}
\end{equation}
To analyze the second term, we consider
\begin{equation*}
\begin{aligned}
& \left(-h^{2}\Delta-\lambda\right)\tilde{\chi}_{2}R\left(\lambda,
h\right)\tilde{\chi}_{1}\\
& \quad =\left(-h^{2}\Delta\tilde{\chi}_{2}\right)R\left(\lambda,
h\right)\tilde{\chi}_{1}-\langle h\nabla\tilde{\chi}_{2}, h\nabla\rangle 
R\left(\lambda,
h\right)\tilde{\chi}_{1}-\tilde{\chi}_{2}\left(h^{2}\Delta-\lambda\right)
R\left(\lambda, h\right)\tilde{\chi}_{1}\\
 & \quad = \left(-h^{2}\Delta\tilde{\chi}_{2}\right)R\left(\lambda,
h\right)\tilde{\chi}_{1}-\langle h\nabla\tilde{\chi}_{2}, h\nabla\rangle 
R\left(\lambda,
h\right)\tilde{\chi}_{1}
\end{aligned}
\end{equation*}
and therefore
\begin{equation}\label{pl2}
\begin{aligned}
& \left\|Q\left(\left(-h^{2}\Delta-\lambda\right)\otimes
I\right)K_{\tilde{\chi}_{2}R\left(\lambda,
h\right)\tilde{\chi}_{1}}\right\|_{L^{2}\left(\mathbb{R}^{2n}\right)}\\
& \quad \leq
C h\left(h\left\|QK_{\left(\Delta\tilde{\chi}_{2}\right)R\left(\lambda,
h\right)\tilde{\chi}_{1}}\right\|_{L^{2}\left(\mathbb{R}^{2n}\right)}+\left\|Q
K_{\langle\nabla\tilde{\chi}_{2}, h\nabla\rangle
R\left(\lambda,
h\right)\tilde{\chi}_{1}}\right\|_{L^{2}\left(\mathbb{R}^{2n}\right)}\right)
\end{aligned}
\end{equation}

Now,
\begin{equation}\label{nablarequiv}
\begin{aligned}
QK_{\langle\nabla\tilde{\chi}_{2}, h\nabla\rangle R(\lambda, 
h)\tilde{\chi}_{1}} &
\equiv\frac{i}{h}\int_{0}^{T}e^{\frac{i}{h}t\lambda}QK_{\langle\nabla\tilde{\chi}_{2},
h\nabla\rangle U(t)\tilde{\chi}_{1}} dt\\
 & = 
\frac{i}{h}\int_{0}^{T}e^{\frac{i}{h}t\lambda}Q(\langle\nabla\tilde{\chi}_{2},
h\nabla\rangle \otimes \tilde{\chi}_{1})K_{U(t)}dt\\
 & = \frac{i}{h}\int_{0}^{T}e^{\frac{i}{h}t\lambda}
\left(Op_{h}(q(\langle\nabla\tilde{\chi}_{2}, \xi\rangle\otimes 
\tilde{\chi}_{1}))
+ h Op_{h}\left(\bar{q}\right)\right) K_{U(t)}dt,\\
\end{aligned}
\end{equation}
where the first equality follows from
(\ref{nowfh}) and $\bar{q}\in C_{c}^{\infty}(\mathbb{R}^{4n}).$
From \[q(\langle\nabla\tilde{\chi}_{2}, \xi\rangle\otimes
\tilde{\chi}_{1})\in S^{0}_{4n}(1)\cap C_{c}^{\infty}(\mathbb{R}^{4n})\]
it follows that $Op_{h}(q(\langle\nabla\tilde{\chi}_{2}, \xi\rangle\otimes
\tilde{\chi}_{1}))=\mathcal{O}_{\mathcal{B}(L^{2}(\mathbb{R}^{2n}))}(1).$
Since $U(t)\in\mathcal{I}_{h}^{0}(\mathbb{R}^{2n}, \Lambda_{t}),$ $t\in\mathbb{R},$ we then have that
\[\left\|Op_{h}(q(\tilde{\chi}_{2}\otimes
\tilde{\chi}_{1})(\tilde{p}\otimes
1))K_{U(t)}\right\|_{L^{2}(\mathbb{R}^{2n})}=\mathcal{O}\left(h^{-\frac{n}{2}}\right)\]
with the norm depending on $t$ continuously.
Therefore, from (\ref{nablarequiv}), we obtain
\begin{equation}\label{nablarestim}
\left\|QK_{\langle \nabla\tilde{\chi}_{2}, h\nabla\rangle R(\lambda,
h)\tilde{\chi}_{1}}\right\|_{L^{2}(\mathbb{R}^{2n})}=\mathcal{O}\left(h^{-\frac{n}{2}-1}\right).
\end{equation}
From
\begin{equation*}
QK_{\left(\Delta\tilde{\chi}_{2}\right)R(\lambda, h)\tilde{\chi}_{1}}
\equiv\frac{i}{h}\int_{0}^{T}e^{\frac{i}{h}t\lambda}QK_{\left(\Delta\tilde{\chi}_{2}\right)U(t)
\tilde{\chi}_{1}}dt
\end{equation*}
we conclude similarly that
\begin{equation}\label{rkere}
\left\|QK_{\left(\Delta\tilde{\chi}_{2}\right)R(\lambda, 
h)\tilde{\chi}_{1}}\right\|
_{L^{2}(\mathbb{R}^{2n})}=\mathcal{O}\left(h^{-\frac{n}{2}-1}\right).
\end{equation}
Also in the same way we obtain
\begin{equation}\label{commutatorkr}
[(P(h)-I)\otimes I, Q]K_{\left(\Delta\tilde{\chi}_{2}\right)R(\lambda,
h)\tilde{\chi}_{1}}=\mathcal{O}_{L^{2}(\mathbb{R}^{2n})}\left(h^{-\frac{n}{2}}\right)
\end{equation}
From (\ref{1pl}), (\ref{pl2}), (\ref{nablarestim}), (\ref{rkere}), 
\eqref{commutatorkr}, and
the fact that
$Op_{h}(g)=\mathcal{O}_{\mathcal{B}(L^{2}(\mathbb{R}^{2n}))}(1)$ for any
$g\in C_{c}^{\infty}(\mathbb{R}^{4n}),$ we then obtain that
\begin{equation}\label{1stpestim}
\left\|Op_{h}(g)(\left(P(h)-\lambda\right) \otimes
I)QK_{\tilde{\chi}_{2}R\left(\lambda,
h\right)\tilde{\chi}_{1}}\right\|_{L^{2}(\mathbb{R}^{2n})}=\mathcal{O}\left(
h^{-\frac{n}{2}}\right),
\end{equation}

Let, now, $f\in C_{c}^{\infty}(\mathbb{R}^{4n})$ also have support near 
$\rho.$
Then
\begin{equation*}
\begin{aligned}
& Op_{h}(f)\left(\left(P(h)-\lambda\right) \otimes I\right)Op_{h}(g)
\left(\left(P(h)-\lambda\right) \otimes 
I\right)QK_{\tilde{\chi}_{2}R\left(\lambda,
h\right)\tilde{\chi}_{1}}\\
& \quad =Op_{h}(f)Op_{h}(g)\left(\left(P(h)-\lambda\right) \otimes
I\right)^{2} QK_{\tilde{\chi}_{2}R\left(\lambda, 
h\right)\tilde{\chi}_{1}}\\
& \quad +Op_{h}(f)[\left(P(h)-\lambda\right) \otimes I,
Op_{h}(g)]\left(\left(P(h)-\lambda\right) \otimes
I\right) QK_{\tilde{\chi}_{2}R\left(\lambda, h\right)\tilde{\chi}_{1}}.
\end{aligned}
\end{equation*}
From (\ref{1stpestim}) and the fact that $[\left(P(h)-\lambda\right) 
\otimes I,
Op_{h}(g)]\in\Psi_{h}^{-1}(1,
\mathbb{R}^{2n})$ we obtain
\begin{equation}\label{2ndpe1}
Op_{h}(f)[\left(P(h)-\lambda\right) \otimes I, Op_{h}(g)]
\left(\left(P(h)-\lambda\right) \otimes I\right)
QK_{\tilde{\chi}_{2}R\left(\lambda,
h\right)\tilde{\chi}_{1}}=\mathcal{O}_{L^{2}(\mathbb{R}^{2n})}\left(h^{1-\frac{n}{2}}\right).
\end{equation}
The same argument as in (\ref{1stpestim}) also implies that
\begin{equation}\label{2ndpe2}
Op_{h}(f)Op_{h}(g)\left(\left(P(h)-\lambda\right) \otimes 
I\right)^{2}QK_{\tilde{\chi}_{2}R\left(\lambda,
h\right)\tilde{\chi}_{1}}=\mathcal{O}_{L^{2}(\mathbb{R}^{2n})}\left(h^{1-\frac{n}{2}}\right)
\end{equation}
and we obtain, from (\ref{2ndpe1}) and (\ref{2ndpe2}),
\begin{equation*}
\left\|Op_{h}(f)\left(\left(P(h)-\lambda\right) \otimes I\right)Op_{h}(g)
\left(\left(P(h)-\lambda\right) \otimes 
I\right)QK_{\tilde{\chi}_{2}R\left(\lambda,
h\right)\tilde{\chi}_{1}}\right\|_{L^{2}(\mathbb{R}^{2n})}=\mathcal{O}\left(h^{1-\frac{n}{2}}\right).
\end{equation*}
Iterating this argument, we then have that
\begin{equation*}
\left\|\left(\prod_{j=j_{2}+1}^{k_{2}}V_{j}\right)QK_{\tilde{\chi}_{2}R\left(\lambda,
h\right)\tilde{\chi}_{1}}\right\|_{L^{2}\left(\mathbb{R}^{2n}\right)}
=\mathcal{O}\left(h^{k-j_{2}-1-\frac{n}{2}}\right), h\to 0.
\end{equation*}

We now observe that
\begin{equation*}
\begin{aligned}
\left(\prod_{j=j_{2}+1}^{k_{2}}V_{j}\right)QK_{\tilde{\chi}_{2}R\left(\lambda,
h\right)\tilde{\chi}_{1}} & =h^{k_{2}-j_{2}}Q K_{\chi_{3}(h)R(\lambda,
h)\tilde{\chi}_{1}}
+h^{k_{2}-j_{2}} K_{\chi_{4}\tilde{P}R(\lambda, h)\tilde{\chi}_{1}}\\
 & \quad
+h^{k_{2}-j_{2}}Op_{h}\left(e_{j_{2}}\right)K_{\tilde{\chi}_{2}\tilde{P}R(\lambda,
h) \tilde{\chi}_{1}},
\end{aligned}
\end{equation*}
where $\chi_{3}(h)$ is polynomial in $h$ with smooth coefficients with 
supports contained in
$\supp \tilde{\chi}_{2},$ $\chi_{4}\in C_{c}^{\infty}(\mathbb{R}^{n})$ is 
such that $\supp
\chi_{4}\subseteq \supp \tilde{\chi}_{2},$ $\tilde{P}\in\Psi_{h}^{0}(1,
\mathbb{R}^{2n})$ with $\sigma\left(\tilde{P}\right)$ is supported in a 
sufficiently small
neighborhood of $\rho,$ $e_{j_{2}}\in C_{c}^{\infty}(\mathbb{R}^{4n}).$
Therefore
\begin{equation*}
\begin{aligned}
&\left(\prod_{i=j_{1}+1}^{j_{2}} V_{i}\right)
\left(\prod_{j=j_{2}+1}^{k_{2}}V_{j}\right)K_{\tilde{\chi}_{2}R\left(\lambda,
h\right)\tilde{\chi}_{1}}\\
 & \quad \equiv
h^{k_{2}-j_{2}-1}i\int_{0}^{T}e^{\frac{i}{h}t\lambda}
\left(\prod_{i=j_{1}+1}^{j_{2}} V_{i}\right) Q 
\left(\chi_{3}(h)\otimes\tilde{\chi}_{1}\right)
K_{U(t)} dt\\
&\quad\quad +
h^{k_{2}-j_{2}-1}i
\int_{0}^{T}e^{\frac{i}{h}t\lambda}\left(\prod_{i=j_{1}+1}^{j_{2}}
V_{i}\right)\left(\chi_{4}\otimes\tilde{\chi}_{1}\right)
(\tilde{P}\otimes I)K_{U(t)} dt\\
&\quad\quad + h^{k_{2}-j_{2}-1} i
\int_{0}^{T}e^{\frac{i}{h}t\lambda}\left(\prod_{i=j_{1}+1}^{j_{2}}
V_{i}\right)Op_{h}\left(e_{j_{2}}\right)\left(\tilde{\chi}_{2}\otimes\tilde{\chi}_{1}\right)
K_{U(t)}dt\\
&\quad 
=\mathcal{O}_{L^{2}(\mathbb{R}^{2n})}\left(h^{k_{2}-1-j_{1}-\frac{n}{2}}\right),
\end{aligned}
\end{equation*}
where we have again used \eqref{nowfh} and the fact that $U(t)\in\mathcal{I}_{h}^{0}(\mathbb{R}^{2n}, \Lambda_{t}),$ $t\in\mathbb{R},$.

Lastly, from the fact that $V_{j}\in\Psi_{h}^{-1}(1, \mathbb{R}^{2n}),$
$1\leq j\leq j_{1},$ we have that
\begin{equation*}
\left(\prod_{j=1}^{k} B_{j}\right)QK_{\tilde{\chi}_{2}R(\lambda, 
h)\tilde{\chi}_{1}}
=h^{k-k_{2}} 
\left(\prod_{j=1}^{k_{2}}V_{j}\right)QK_{\tilde{\chi}_{2}R(\lambda,
h)\tilde{\chi}_{1}}=\mathcal{O}_{L^{2}(\mathbb{R}^{2n})}\left(h^{k-1-\frac{n}{2}}\right),
\end{equation*}
which completes the proof of the theorem.  $\hfill\Box$

\subsection{Non-Trapping Energy Level}\label{snt}

\begin{Th}\label{ntfio}
Let $\lambda>0$ be a non-trapping energy level for $P.$

Then
amplitude $A\left(\lambda, h\right)\in
\mathcal{I}_{h}^{\frac{3}{2}}(\mathbb{S}^{n-1}\times\mathbb{S}^{n-1},
SR \left(\lambda\right)).$
\end{Th}

\medskip
\noindent
{\it Proof:}\;
We recall from \cite{VZ} that $\|\tilde{\chi}_{2}
R\left(\lambda, h\right)
\tilde{\chi}_{1}\|_{\mathcal{B}(L^{2}(\mathbb{R}^{n}))}=\mathcal{O}\left(\frac{1}{h}\right).$
Then the result follows from Theorem \ref{resolventfio} and the main
theorem. $\hfill\Box$

\subsubsection{A Simple Inverse Problem}\label{sinv}

Following a suggestion of Plamen Stefanov we include a discussion of
an inverse problem motivated by Theorem \ref{ntfio}. Suppose that
$ P = - h^2 \Delta + V $, where $ V $ is compactly supported and smooth,
that satisfies the general assumptions of this article.
Suppose that
$ \lambda > \max_{x \in {\mathbb R}^n } V_+ ( x ) $ so that the energy
level $ \lambda $ is clearly non-trapping.
Let $V$ further be such that the metric $(\lambda-V(x))dx^2$ is simple.
We have the following

\begin{Th}\label{t:stefanov}
For $ P$ and $ \lambda $ as above, the scattering relation, and hence by
Theorem \ref{ntfio}
the scattering amplitude, determine $ V $ uniquely.
\end{Th}

\medskip
\noindent
{\em Outline of proof:}
We compare
this with the problem of the wave equation with variable wave speed 
defined as
\[ c(x) = ( \lambda - V( x))^{-\frac12} \,.\]
Then the function $c$ is equal to $ 1 / \sqrt{\lambda} $  for large $x.$
We have a new Hamiltonian $\widetilde p -1
=c^{-2}(x)|\xi|^2 -1$. These Hamiltonians have the same integral
curves but they are parameterized in different ways. Hence
the scattering relation
for this Schr\"odinger equation is that related to the metric 
$g_1(x)=c^{-2}(x)dx^2.$
It is now implicit in \cite{MichelR} that the scattering relation for the 
metric $g_1$
determines the boundary distance function uniquely.
The results of \cite{MuR} and \cite{Mu} further imply that the boundary 
distance function
determines uniquely a simple metric conformal to the Euclidean, in 
particular, it determines
$c$ and therefore $V.$
$\hfill\Box$

\subsection{Trapping Energy Level}\label{str}
Here we consider a trapping energy level $\lambda>0.$
We make the following assumption
\begin{Assume}\label{resassumption}
\begin{equation*}
\begin{aligned}
& \exists\; \epsilon>0 :\text{ the resonances } \lambda_{j} \text{ satisfy 
}\\
& |\Im \left(\lambda_{j}\right)|\geq Ch^{q} \text{ for } \Re
\left(\lambda_{j}\right)\in \left[\lambda-\epsilon,
\lambda+\epsilon\right].
\end{aligned}
\end{equation*}
\end{Assume}

We let $V\subset\Lambda_{R}(\lambda)$ be as in Theorem \ref{resolventfio}.
Then for every $\left(x, \xi\right)\in \pi_{1}(V)$ there exist unique 
$t_j\left(x,
\xi\right)\in\mathbb{R}$ such that
$\exp\left(t_j\left(x, \xi\right)H_p\right)\left(x, \xi\right)\in
L_j\left(\lambda\right),$ $j=1, 2.$
As in Section \ref{scr} we have that there exists an open set
$U\subset\cup_{\left(x, \xi\right)\in \pi_{1}(V)}\exp\left(t_1\left(x,
\xi\right)H_p\right)\left(x, \xi\right)$ such that we
can define the scattering relation $SR_U\left(\lambda\right)$ as in
Section
\ref{scr}.
By decreasing $U,$ if necessary, we can further assume that
$SR_{\bar{U}}\left(\lambda\right)$ is a Lagrangian submanifold of
$T^{*}\mathbb{S}^{n-1}\times T^{*}\mathbb{S}^{n-1}.$

Under these assumptions we have the following
\begin{Th}\label{tfio}
$A(\lambda, h)\in
\mathcal{I}_{h}^{\frac{3}{2}}(\mathbb{S}^{n-1}\times\mathbb{S}^{n-1},
SR_{\bar{U}}\left(\lambda\right)).$
\end{Th}

\medskip
\noindent
{\it Proof:}\;
Let $C\in\Psi_{h}^{0}\left(1,
\mathbb{S}^{n-1}\times\mathbb{S}^{n-1}\right)$ have
compactly supported
symbol which does not vanish in a neighborhood of a point $\rho\in
SR_{\bar{U}}\left(\lambda\right).$

We recall the following result contained in the proof of Proposition 5.1,
\cite{Michel}.

\begin{Prop}
Let Assumption \ref{resassumption} hold.
Then there exists $\tilde{n}\in \mathbb{N}$ such that
\begin{equation*}
\|\chi R\left(z,
h\right)\chi\|_{\mathcal{B}(L^{2}(\mathbb{R}^{n}))}=\mathcal{O}\left(h^{-\tilde{n}}\right),
h\to 0,
\end{equation*}
uniformly for $z\in\left\{z\in\mathbb{C}: \left(\Re z, \pm\Im z\right)\in
\left(\frac{1}{d}, d\right)\times(0,1]\right\}$ and $\chi\in
C_{c}^{\infty}\left(\mathbb{R}^{n}\right),$ $\chi=1$ on $B\left(0,
\rho\right),$ for $\rho$
large enough.
\end{Prop}

This Proposition was proved in \cite{Michel} for a certain class of
potential
perturbations
of the Euclidean Laplacian.
The proof there, however, relies only on Lemma 4.1, \cite{TZ}, and the
latter is proved in the black-box setting.
Therefore the above Proposition holds also in the setting discussed here.

From Theorem \ref{resolventfio} we have that $\tilde{\chi}_{2}R(\lambda,
h)\tilde{\chi}_{1}\in\mathcal{I}_{h}^{1}\left(\mathbb{R}^{2n},
\overline{V}\cap\Lambda_{R}(\lambda)\right)$ and the assertion of the
Theorem now follows by the main theorem.$\hfill\Box$

\subsection{Microlocal Representation of the Scattering
Amplitude}\label{smicrol}

Here we prove the following

\begin{Th}\label{tmicrol}
Let $\omega_0\in\mathbb{S}^{n-1}$ be regular for 
$\theta_0\in\mathbb{S}^{n-1}$ and
$L\in\mathbb{N}$ be the number of $(\theta_0, \omega_0)$ phase 
trajectories.
Let $P_l\in\Psi_{h}^{0}(1, \mathbb{S}^{n-1}\times\mathbb{S}^{n-1}),$ $l=1,
\dots, L,$ be microlocal cut-offs to the Lagrangian submanifolds 
$SR_l(\lambda)$ defined by
\eqref{lagrl},
respectively.

Then
\begin{equation*}
P_l A(\lambda, h)=e^{\frac{i}{h}S_{l}}a_{l},\: l=1, \dots, L,
\end{equation*}
where $S_{l},$ $l=1, \dots, L,$ are as given by
(\ref{maction}) and $a_l\in S_{2n-2}^{\frac{n}{2}+1}(1),$ $l=1, \dots, L,$
have compact support.
\end{Th}

\medskip
\noindent
{\it Proof:}\;
By Theorem \ref{ntfio} and Theorem \ref{tfio}, the
scattering
amplitude
is a global semi-classical Fourier integral operator associated to
$\cup_{l=1}^{L}SR_{l}(\lambda).$
The assertion of the Theorem then follows from Lemma \ref{actionsr}, \cite[Lemma 5]{AI},
and \cite[Theorem 1]{AI}. $\hfill\Box$

We remark here that the phase function in this microlocal representation 
of the scattering
amplitude is the same as the one given by \cite{Michel} and \cite{RT}.

{\bf Acknowledgements.}
I would like to thank Maciej Zworski for supervising my Ph.~D. thesis of 
which this paper formed a part.
I would also like to thank Vesselin Petkov for introducing me to the 
Ph.~D. thesis  of his student Laurent Michel, which has helped me complete my work on this project and Plamen Stefanov for discussions leading to the material presented in Section \ref{sinv}.
I am further grateful to Victor Ivrii and Xiang Tang for helpful discussions.

\end{document}